\documentclass[11pt,a4paper]{article}
\usepackage{ifthen}
\usepackage{amsmath}
\usepackage{amssymb}
\usepackage{latexsym}
\usepackage{bbm}
\usepackage{graphicx}
\usepackage{pdfsync}
\usepackage{verbatim}
\usepackage{fancyhdr}
\usepackage{url}               
\usepackage[english]{babel}
\usepackage[applemac]{inputenc}

\usepackage[pdftex,%
a4paper=true,%
bookmarks=true,%
bookmarksnumbered=true,%
colorlinks=true,%
urlcolor=blue,%
pdfstartview=FitH%
]{hyperref}

\hypersetup{
	colorlinks=true,     
	citecolor=blue,
	linkcolor=red,     
	}

\usepackage{ntheorem}
\newtheorem{theorem}{Theorem}[section]
\newtheorem{lemma}[theorem]{Lemma}

\newtheorem{proposition}[theorem]{Proposition}
\newtheorem{corollary}[theorem]{Corollary}
\newtheorem{predefinition}[theorem]{Definition}
\newenvironment{definition}{\begin{predefinition}\upshape}{\end{predefinition}}
\newtheorem{preremark}[theorem]{Remark}
\newenvironment{remark}{\begin{preremark}\upshape}{\end{preremark}}

\numberwithin{equation}{section}
\numberwithin{figure}{section}

\newcommand{\bp}{{\it Proof. }}
\newcommand{\ep}{\hfill $\square$\\}

\newcommand{\bpl}{{\it Proof of Lemma }} 
\newcommand{\epl}{\hfill $\square$\\}

\newcommand{\bpp}{{\it Proof of Proposition }} 
\newcommand{\epp}{\hfill $\square$\\}

\newcommand{\bpc}{{\it Proof of Corollary }} 
\newcommand{\epc}{\hfill $\square$\\}

\newcommand{\Imm}{\mbox{Im}}
\renewcommand{\Re}{\mbox{Re}}
\newcommand{\be}{\begin{equation}}
\newcommand{\ee}{\end{equation}}
\newcommand{\bea}{\begin{eqnarray}}
\newcommand{\eea}{\end{eqnarray}}
\newcommand{\bee}{\begin{eqnarray*}}
\newcommand{\eee}{\end{eqnarray*}}

\def\pa{\partial}

\def\na{\nabla}

\def\RR{\mathbb{R}}
\def\ZZ{\mathbb{Z}}

\def\ni{\noindent}

\def\eps{\vare}

\def\eps{\varepsilon}

\catcode`@=11
\def\supess{\mathop{\operator@font Sup\,ess}}
\catcode`@=12

\def\RR{\mathbb{R}}

\def\ZZ{\mathbb{Z}}

\def\ds{\displaystyle}

\def\ni{\noindent}

\def\bar#1{{\overline #1}}

\def\R2+{\RR ^2_+}

\def\pa{\partial}
\def\na{\nabla}

\def\lim{\mathop{\rm lim}}

\def\supp{{\rm supp}~}
\def\sup{\mathop{\rm sup}}
\def\exp{{\rm exp}}

\def\pa{\partial}

\newcommand{\uD}{\mathrm{D}}
\newcommand{\ud}{\mathrm{d}}
\newcommand{\ue}{\mathrm{e}}
\newcommand{\ui}{\mathrm{i}}

\begin{document}

\renewcommand{\refname}{References}
\bibliographystyle{alpha}

\pagestyle{fancy}
\fancyhead[L]{ }
\fancyhead[R]{}
\fancyfoot[C]{}
\fancyfoot[L]{ }
\fancyfoot[R]{}
\renewcommand{\headrulewidth}{0pt} 
\renewcommand{\footrulewidth}{0pt}

\newcommand{\montitre}{Gevrey regularity and analyticity for the solutions of the Vlasov-Navier-Stokes system}

\newcommand{\auteur}{\textsc{ Dahmane Dechicha
}}
\newcommand{\affiliation}{Laboratoire J.-A. Dieudonn\'e. Universit\'e C\^{o}te d'Azur\\ Parc Valrose, 06108 Nice Cedex 02,   France\\
Laboratoire de recherche AGM. CY Cergy Paris Universit\'e. UMR CNRS 8088 \\
2 Avenue Adolphe Chauvin 95302 Cergy-Pontoise Cedex, France \\
\url{dahmane.dechicha@cyu.fr}}

 \begin{center}
{\bf  {\LARGE \montitre}}\\ \bigskip \bigskip
 {\large\auteur}\\ \bigskip \smallskip
 \affiliation \\ \bigskip
\today
 \end{center}
 \begin{abstract}
In this paper, we prove propagation of $\frac{1}{s}$-Gevrey regularity $(s \in (0, 1))$ and analyticity $(s=1)$ for the Vlasov-Navier-Stokes system on  $\mathbb{T}^d \times \RR^d$ (and $\RR^d\times\RR^d$) using a Fourier space method in analogy to the results proved for the Euler system in \cite{KV1} and \cite{LO} and for Vlasov-Poisson system in \cite{VR}. More precisely, we give quantitative estimates for the growth of the $\frac{1}{s}$-Gevrey norm and decay of the regularity radius for the solution of the system in terms of $\na_x u$,  the spatial  density $\rho_f$ and the diameter of the support in the velocity variable of the distribution of particles $f$. 
 In particular, this implies existence of $\frac{1}{s}$-Gevrey $(s \in (0, 1))$ and analytic $(s = 1)$ solutions for the Vlasov-Navier-Stokes system in $\mathbb{T}^d\times\RR^d$ (and $\RR^d\times\RR^d$),  and global Gevrey solutions in $\mathbb{T}^3\times\RR^3$  for sufficiently small data, and an initial data for the Vlasov equation with compact support in velocity.
\end{abstract}





\tableofcontents


 \pagestyle{fancy}
\fancyhead[R]{\thepage}
\fancyfoot[C]{}
\fancyfoot[L]{}
\fancyfoot[R]{}
\renewcommand{\headrulewidth}{0.2pt} 
\renewcommand{\footrulewidth}{0pt} 

\thispagestyle{empty}

\section{Introduction}
\subsection{Setting of the problem and historical context}
In this paper, we deal with the Vlasov-Navier-Stokes system (VNS):
\begin{equation}\label{VNS}
(\mathrm{VNS})\left\{\begin{array}{rcl}
\partial_t f + v\cdot\nabla_x f + \nabla_v\cdot[(u-v)f] &=& 0 , \hspace{1.95cm} \mbox{in } [0,T)\times\mathbb{T}^d\times\mathbb{R}^d, \\
\partial_t u + (u\cdot\nabla_x)u - \Delta_x u + \nabla_x p &=& j_f - \rho_f u , \hspace{0.6cm}  \mbox{ in } [0,T)\times\mathbb{T}^d,\\
\nabla_x\cdot u &=&0 , \hspace{1.85cm} \mbox{ in } [0,T)\times\mathbb{T}^d, \\
u(0,\cdot)= u_0, \quad f(0,\cdot,\cdot) &=& f_0 ,
\end{array}\right.
\end{equation}
where $T\in]0,\infty]$, and where $\rho_f$ and $j_f$ are the spatial density and local current respectively:
$$  \rho := \rho_f(t,x) := \int_{\mathbb{R}^d} f \ \ud v \quad  \mbox{ and } \quad  j :=j_f(t,x) := \int_{\mathbb{R}^d}v f \ \ud v .  $$

This system of nonlinear PDEs describes the transport of particles (described by their density function $f$) within a homogeneous incompressible fluid (described by its velocity $u$ and its pressure $p$).  This description corresponds to a regime where the particles volume fraction is small compared to that of the surrounding fluid.  It belongs to the broad family of \emph{fluid-kinetic systems} or \emph{couplings},  which were introduced in the pioneering works of O'Rourke \cite{o1981collective} and Williams \cite{williams1985combustion} for the description of sprays involving a large number of particles.  We also refer to \cite{desvillettes2010some} for a general overview on the description of multiphase flows, as well as to \cite{reitz1996computer}.  Among all possible couplings (we refer to the introduction of \cite{glass2018vlasov} or \cite{ertzbischoff2023well} for other examples), the Vlasov-Navier-Stokes system has been intensively studied because of both its physical relevance and the mathematical challenges that it offers.  It has been for instance shown to provide a good description of medical aerosols in the upper part of the lung (see e.g. \cite{boudin2015modelling,BM21}). The VNS system \eqref{VNS} is fully coupled: both unknowns $f$ and $u$ depend on each other. This is due to the Brinkman force (the source term in the fluid equation) and the drag acceleration (the inertial term in the kinetic equation). We refer to \cite{boudin2015modelling} for the physical justification of these, and to \cite{desvillettes2008mean,BDGR1, BDGR2, hillairet2018homogenization, hillairet2019effect}  for the (partial) mathematical derivation of the former. The physical constants are all normalized in \eqref{VNS}. The VNS system can also be considered with inhomogeneous or compressible Navier-Stokes equations \cite{choi2015global, choi2017finite, choi2022regular} and additional terms in kinetic equations \cite{choi2022temporal, choi2022local}.
Note that the case of compressible Euler equations for the fluid, coupled to a kinetic equation, has also been investigated  \cite{baranger2006coupling}.   

The study of the VNS system, from a mathematical point of view, has been the topic of several research papers in the last twenty years and in many directions of research.  The Cauchy theory, addressing the existence of weak global solutions for the VNS system, has been tackled in dimension 2 and 3 in various domains of space (see for instance \cite{anoshchenko1997existence,hamdache1998global,BDGM}), and also allows for more complex physics in the model (see \cite{boudin2017global,boudin2020global}).  It mainly consists in obtaining a \emph{Leray weak} solution for $u$ and a \emph{renormalized weak} solution (in the sense of Di-Perna and Lions \cite{diperna1989ordinary}) for $f$, using a remarkable energy-dissipation identity that is satisfied by solutions to the system. In dimension 2, the uniqueness of such solutions has been shown in \cite{han2019uniqueness}.  In \cite{choi2015global}, Choi and Kwon showed the existence of a unique strong solutions to the inhomogeneous VNS system in a time interval that depends on the initial data (provided that the initial data is sufficiently small and regular),  and also established an a priori estimate for the large-time behavior of the solutions to the last system in a spatial periodic domain, i.e.  in $\mathbb{T}^3$.   This last question, concerning the long-time behavior of (VNS) solutions and large-scale dynamics, has attracted a lot of attention and has been the subject of several advances over the last few years.  Indeed, it is expected that the cloud of particles aligns its velocity on that of the fluid,  
$$ u(t) \underset{t \to +\infty}{\longrightarrow} v^\infty , \quad f(t) \underset{t \to +\infty}{\longrightarrow} \rho^\infty \otimes \delta_{v=v^\infty} ,  \footnote{ \mbox{ The convergence of} $u(t)$ \mbox{is in} $L^2$, \mbox{while that of} $f(t)$ \mbox{is in the sense of the Wasserstein(-1) distance} $W_1$.} $$
for some asymptotic velocity $v^\infty \in \RR^3$ and profile $\rho^\infty(t, x)$.  The answer to this question was obtained first in \cite{HMM} for Fujita-Kato type solutions, in the $3d$ torus case and in \cite{han2022large} in the whole space $\RR^3$, while the case of a $3d$ bounded domain (with absorption boundary conditions for $f$) and the half-space case are investigated in \cite{ertzbischoff2021concentration} and \cite{ertzbischoff2021decay} respectively.  Another type of asymptotics has been studied on VNS, we refer to \cite{han2021hydrodynamic,hofer2018inertialess,ertzbischoff2022global} for more details.   It should be noted that the ``rigorous'' derivation of the VNS system from an ODE system written at the microscopic scale remains an open problem. We refer to \cite{desvillettes2008mean, hillairet2018homogenization, hillairet2019effect, carrapatoso2020derivation, hillairet2021derivation, hofer2020convergence} for a partial answer.    

In this work, we investigate the propagation of higher regularity for smooth solutions to the Vlasov-Navier-Stokes system.  Our method is based on the notion of \emph{Gevrey class} regularity, which is a stronger concept than the $C^\infty$ regularity.  It not only asserts that all derivatives of the solution are bounded, but also that these bounds depend on the order of the derivatives in some prescribed way. Gevrey \cite{gevrey1918nature} used this notion as a setting in which to extend Cauchy-Kowalevski existence arguments to classes of functions that are not necessarily analytic (for a review of the analytic case see, e.g., \cite{john1982partial}). In fact, they are special cases of the quasianalytic classes \cite{hadamard1912generalisation}.  La Vall\'ee Poussin \cite{la1924quatre} showed that, among the quasianalytic functions, the Gevrey classes are characterized by an exponential decay of their Fourier coefficients,  see \cite[Section 2]{LO} for a definition of Gevrey classes and the proof of their characterization by Fourier transformation and Sobolev spaces (see also \cite{kopec1960quasianalytic}).  An equivalent definition is given in Subsection \ref{notations et definitions}.  In turn, this characterization has proven useful in investigating different questions for the solutions of various nonlinear partial differential equations.  For example, this notion played a very important role in the proof of \emph{Landau damping} in the paper by Bedrossian, Masmoudi and Mouhot \cite{BMM}, where the authors improved the Mouhot-Villani result \cite{MV}.    

Concerning the question of propagation of higher regularity,  more precisely the propagation of Gevrey regularity and analyticity, has been the subject of several papers since the seventies. In \cite{FoiTem}, Foias and Temam proved that the solution of the Navier-Stokes equations (in dimension 2 and 3) is in the Gevrey class  for a Sobolev initial data and a Gevrey source term which does not depend on the solution, with an affine (in time) radius of regularity. The fact that the source term here is not complicated (since it is independent of the solution) plays an important role in this gain of regularity and it allows to take advantage of the dissipative term.  Contrary to this last result, for Euler's equations, the radius of analyticity (of regularity in the Gevrey case) decays exponentially in time (and as long as the quantity $\int_0^t \| \na u(s)\|_ \infty \ud s$, which appears in the exponent, remains finite).  The first result on the Euler system have been derived by Bardos, Benachour and Zerner (see references \cite{BBZ,Bardos1976,Benachour}) who use estimates on the Green function of the Poisson kernel in the complex plane to describe the region of analyticity.  For results on the local propagation of analyticity see the works by Baouendi and Goulaouic \cite{baouendi1975problemes}, Alinhac and Metivier \cite{alinhac1985propagation, alinhac1986propagation}, and references therein.   These results were continued by the work of Levermore and Oliver \cite{LO} where they proved the propagation of the analytic regularity in dimension 2 on $\mathbb{T}^2$ using a method of Fourier space based on the notion of Gevrey regularity. However, the analyticity radius decay rate obtained by these last two authors was $\exp(-\exp(\exp t))$, which is faster than the $\exp(-\exp t)$ obtained by Bardos, Benachour and Zerner in \cite{BBZ}. Subsequently, Levermore and Oliver's result was improved by Kukavica and Vicol in \cite{KV1} who showed the same rate of decay for the analyticity radius obtained in \cite{BBZ}, but using a Fourier space method instead.   

Recently, using a Fourier space method in analogy to the results proved for the $2d$-Euler system in \cite{KV1} and \cite{LO} and applying techniques used in the proof of Landau damping \cite{BMM}, Velozo Ruiz \cite{VR} proved the Gevrey regularity propagation for solutions of the Vlasov-Poisson system,  giving a quantitative estimate of the decay in the radius of regularity (it is an  $\exp(-\exp t)$ decrease which was obtained, as for Euler) for the solution of the system in terms of the force field and the volume of the support in the velocity variable of the distribution of matter.    

In this paper we address the problem of propagation of Gevrey regularity for the VNS system on $\mathbb{T}^d \times \RR^d$ as long as there exists a Sobolev solution $(f,u)$ for this system. More precisely, we give quantitative estimates for the growth of the Gevrey norm and decay of regularity radius for the solution of \eqref{VNS} in terms of the Sobolev norm which is itself estimated in terms of $\| u\|_{W^{1,\infty}}$, the density $\| \rho_f\|_\infty$ and the diameter of the support in the velocity variable of the distribution $f$.  As an application, we show global existence of Gevrey solutions for the VNS system in $\mathbb{T}^3\times \RR^3$ for initial \emph{modulated energy} small enough, due to the result proved by Han-Kwan,  Moussa and Moyano in \cite{HMM}. Furthermore,  the propagation of Gevrey regularity remains true on $\RR^d\times\RR^d$ even if it means replacing the Fourier series by integrals.

\subsection{Notations,  definitions and preliminaries}\label{notations et definitions}
In order to write the main theorems of the paper,  let introduce the usual Gevrey norms. \\
In the following, we use the multi-index notations 
$$ v^\alpha := (v_1)^{\alpha_1}\ldots(v_d)^{\alpha_d} \quad \mbox{ and }  \quad \uD^\alpha_\eta := (\ui\partial_1)^{\alpha_1} \ldots(\ui\partial_d)^{\alpha_d} ,$$
where $\alpha = (\alpha_1,\ldots,\alpha_d) \in \mathbb{N}^d, v \in \mathbb{R}^d$, $\eta\in \mathbb{R}^d$ and $\ \ui^2=-1$. \\
We define the usual Fourier coefficient (transformation) of $f \in L^2(\mathbb{T}^d\times\RR^d)$ by
$$ \mathcal{F}(f)(k,\eta) :=   \hat f_k(\eta) := \frac{1}{(2\pi)^d} \iint_{\mathbb{T}^d\times\RR^d} \ue^{-\ui x\cdot k-\ui v\cdot \eta} f(x,v) \ud x \ud v $$
and of $u \in L^2(\mathbb{T}^d)$ by
$$  \mathcal{F}(u)(k) :=   \hat u_k := \frac{1}{(2\pi)^d} \int_{\mathbb{T}^d} \ue^{-\ui x\cdot k} u(x) \ud x .$$
We denote by $\langle \cdot,\cdot \rangle_{L^2}$ the scalar product in the Hilbert space $L^2$ and
we define the Japanese brackets: $\ds \langle k \rangle := (1+|k|^2)^{\frac{1}{2}}$ and $\ds \langle k,\eta \rangle := (1+|k|^2+|\eta|^2)^{\frac{1}{2}}$ for all $k, \eta \in \RR^d$. \\
Finally, we denote the standard Sobolev norm of $f$ in $H^\sigma_{x,v}(\mathbb{T}^d\times\RR^d)$ by $\|f\|_\sigma$ and we denote by  $H^\sigma_{x,v;M}(\mathbb{T}^d\times\RR^d)$ (as in \cite{VR}) the weighted Sobolev space with the norm 
$$ \| f\|_{\sigma,M}^2 := \sum_{|\alpha|\leqslant M} \| v^\alpha f \|_{\sigma}^2 ,$$
which can be written, in Fourier variables, as 
$$ \| f\|_{\sigma,M}^2 := \sum_{|\alpha|\leqslant M} \sum_{k \in \ZZ^d}  \int_{\RR^d} \big|\uD^\alpha_\eta \hat f_k(\eta) \big|^2 \langle k,\eta\rangle^{2\sigma} \ud \eta .$$
\begin{definition}\textbf{($\frac{1}{s}$-Gevrey Classes in $\mathbb{T}^d\times\mathbb{R}^d$).}
A real-valued function $f \in C^\infty(\mathbb{T}^d\times\mathbb{R}^d)$ is said to be of Gevrey class $\frac{1}{s}$ with radius of regularity $\lambda>0$, Sobolev correction $\sigma>0$ and weight $M\in \mathbb{N}$, if for some $s\in (0,1]$, we have $f \in L^2(\mathbb{T}^d\times\mathbb{R}^d)$ and
$$\| f \|^2_{\lambda,\sigma,M,s} := \underset{|\alpha|\leqslant M}{\sum} \| v^\alpha f \|^2_{\lambda,\sigma,s}< +\infty ,$$
with
$$\| v^\alpha f \|^2_{\lambda,\sigma,s} := \| A\uD^\alpha_\eta \hat{f} \|^2_{L^2_{k,\eta}} := \underset{k\in\mathbb{Z}^d}{\sum}\int_{\RR^d}\langle k,\eta\rangle^{2\sigma}\ue^{2\lambda\langle k,\eta\rangle^s}|\uD^\alpha_\eta \hat{ f}_k(\eta)|^2 \ud \eta, $$
and where $$ A:= A^\sigma_k(\eta)=\langle k,\eta\rangle^{\sigma}\ue^{\lambda\langle k,\eta\rangle^s}$$
is the Fourier multiplier. We denote by $\mathcal{G}^{\lambda,\sigma,M,\frac{1}{s}}(\mathbb{T}^d\times\RR^d)$ the space of functions of this class.
\begin{remark} 
The standard definition of Gevrey spaces does not include the term  $v^\alpha$.  This term is equivalent to $(1+|v|^M)$, and plays the role of a weight for the velocity variable. It was introduced to control the Sobolev/Gevrey norms of $\rho_f$ and $j_f$ by those of $f$. 
\end{remark}
\end{definition}
\begin{definition}\textbf{($\frac{1}{s}$-Gevrey Classes in $\mathbb{T}^d$).}
A real vector function $u\in C^\infty(\mathbb{T}^d;\mathbb{R}^d)$ is said to be of Gevrey class $\frac{1}{s}$ with radius of regularity $\lambda>0$ and Sobolev correction $\sigma>0$ if,  for some $s\in (0,1]$, we have $u \in L^2(\mathbb{T}^d)$ and
$$\| u \|^2_{\lambda,\sigma,s} := \| \ue^{\lambda\Lambda^s}u \|^2_{\sigma} :=\| \Lambda^\sigma \ue^{\lambda\Lambda^s}u \|^2_2 := \underset{k\in\mathbb{Z}^d}{\sum}\langle k\rangle^{2\sigma}\ue^{2\lambda\langle k\rangle^s}|\hat{u}_k|^2 < +\infty , $$
where $$ \Lambda := (\mathrm{Id} - \Delta_x)^{\frac{1}{2}}.$$
In Fourier variables, $\displaystyle  \langle k\rangle^{\sigma }e^{\lambda\langle k\rangle^s}=:A^\sigma_k(0)$ is the Fourier multiplier and $\hat{u}_k$ the Fourier coefficients of $u$ on $\mathbb{T}^d$. We denote by $\mathcal{G}^{\lambda,\sigma,\frac{1}{s}}(\mathbb{T}^d)$ the space of functions of this class.
\end{definition}

\ni \textbf{The transport equation.} The Vlasov equation
$$ \pa_t f + v \cdot \na_x f + \na_v \cdot \big[ (u-v) f \big] = 0 $$
can be rewritten as
$$  \pa_t f + v \cdot \na_x f + (u-v) \cdot \na_v f  -d f = 0  . $$

For $s, t \geqslant 0$ and $(x, v) \in \mathbb{T}^d\times\RR^d$, we define (see \cite[Definition 4.1]{HMM}) the characteristic curves $(X(s, t, x, v), V (s, t, x, v))$ as the solutions to the system of ODEs  
\begin{equation}\label{caracteristiques}
\left\{\begin{array}{rcl} 
\ds \frac{\ud}{\ud s} X(s,t,x,v) &=& V(s,t,x,v) ,   \hspace{4.329cm} X(t,t,x,v)=x , \\
\\
\ds \frac{\ud}{\ud s} V(s,t,x,v) &=& u(s,X(s,t,x,v)) - V(s,t,x,v) ,  \hspace{1.1cm} V(t,t,x,v)=v .
\end{array}\right.
\end{equation}
By the method of characteristics, for a smooth vector field $u$, we can write the solution $f$ to the Vlasov equation as
\begin{equation}\label{f = e^(dt)f_0}
f(t,x,v) = \ue^{dt} f_0(X(0,t,x,v),V(0,t,x,v)) .
\end{equation}
As a consequence,  for almost all $t \geqslant 0$,
$$ \|f(t) \|_{L^\infty(\mathbb{T}^d\times\RR^d)} \leqslant  \ue^{dt} \|f_0 \|_{L^\infty(\mathbb{T}^d\times\RR^d)} .$$
For $f_0$ with compact support, thanks to \eqref{f = e^(dt)f_0} and under certain conditions on the regularity and smallness of $u$, the functions $X$ and $V$ behave well and therefore, the function $f(t,\cdot,\cdot)$ will also have a compact support for any $t$. We refer to Section 4 in \cite{HMM} for more details, and to \cite[Lemma 2.3]{choi2015global} for estimating the support of $f$ for small data. In this paper, we do not need to study equations \eqref{caracteristiques} but only have an estimate of the quantity $\| f (t) \|_{\infty,M}$. In Lemma \ref{lemma f_inf,M}, we give a control of the last norm in terms of the diameter of the support of $f$ and $\|f_0\|_{\infty,M}$.  Specifically, for $u \in L^1(0,t;L^\infty)$, one has: $\| f (t) \|_{\infty,M} \lesssim \ue^{dt} \big(1+\|u\|_{L^1(0,t;L^\infty)}^M\big)\|f_0\|_{\infty,M}$,  and for $d=3$, this simplifies to: $\| f (t) \|_{\infty,M} \lesssim \ue^{dt} \|f_0\|_{\infty,M}$.  \\  

\ni \textbf{On the existence of strong solutions for VNS. }
In \cite{choi2015global}, the existence of a unique strong solution was proved for the inhomogeneous Vlasov-Navier-Stokes system in $\Omega\times\RR^3$ with $\Omega=\mathbb{T}^3$ or $\RR^3$,  under some assumptions on the density $\varrho$,  taking $f_0$ with compact support in position and velocity and under the smallness of $\|f_0\|_{H^2(\Omega\times\RR^3)} + \|u_0\|_ {H^2(\Omega)}$. For our system \eqref{VNS} which corresponds to $\varrho \equiv 1$ in \cite[Theorem 1.1]{choi2015global}, the solution is given in the following sense:\\  
For any $T >0$, there exists $\eps:=\eps(T) >0$ depending only on $T$ such that if
$$ \|f_0\|_{H^2(\Omega\times\RR^3)} + \|u_0\|_ {H^2(\Omega)} \leqslant \eps ,  $$
then, the VNS system \eqref{VNS} admits a unique strong solution $(f,u)$ satisfying
\begin{itemize}
\item $f \in C(0,T;H^2(\Omega\times\RR^3))$ ;
\item $u \in C(0,T;H^2(\Omega)) \cap L^2(0,T;H^3(\Omega)) \ $ and $ \ \pa_t u \in C(0,T;L^2(\Omega)) \cap L^2(0,T;H^1(\Omega))$ ;
\item $\na_x p \in C(0,T;L^2(\Omega)) \cap L^2(0,T;H^1(\Omega))$.
\end{itemize}

A \emph{modulated} version, denoted by $\mathcal{E}(t)$, of the energy $E(t)$ (see Definition \ref{energie modulee} and Definition \ref{energie} respectively) was introduced in the paper \cite{choi2015global} and played an important role in the work of \cite{HMM}.  In particular,  for $\mathcal{E}(0)$ defined by
\begin{align}\label{energie modulee}
\mathcal{E}(0):= \frac{1}{2} \int_{\mathbb{T}^3\times\RR^3} f_0(x,v)&|v-\langle j_{f_0}\rangle|^2 \ \ud v \ud x +\frac{1}{2} \int_{\mathbb{T}^3} |u_0(x)-\langle u_0\rangle|^2 \  \ud x + \frac{1}{4} \big|\langle j_{f_0}\rangle- \langle u_0\rangle\big|^2 ,
\end{align}
small enough,  in the sens of \eqref{petitesse de E}, we get
\begin{equation}\label{consequence de E(0)}
\int_1^\infty \| \na_x u(t) \|_\infty \ud t \ll 1  \ \mbox{ and } \  \| \rho_f \|_{L^\infty((0,\infty)\times\mathbb{T}^3)} \lesssim 1 ,
\end{equation}
where $\langle j_{f_0}\rangle := \int_{\RR^3} j_{f_0}(x) \ud x$ and $\langle u_0\rangle := \int_{\RR^3} u_0(x) \ud x$. 
This last estimate will allow us to obtain the global existence of Gevrey solutions to the VNS system \eqref{VNS} in $\mathbb{T}^3\times\RR^3$ later on.

\subsection{Main results}\label{main results}
The results of this paper hold in any dimension $d$, except for Corollary \ref{existence globale} which is for $d=3$.
From now on, the parameter $s \in (0,1]$ is fixed, while $\lambda(t)$ can vary over time.

\begin{theorem}[Propagation of Gevrey regularity]\label{propagation Gevrey}  Let $(f_0,u_0)$ be initial data for the VNS system \eqref{VNS} on $\mathbb{T}^d\times\RR^d$ (or $\RR^d\times\RR^d$) such that,  $f_0$ has a compact support in velocity and  $\|f_0\|_{\lambda_0,\sigma,M,s} + \|u_0\|_{\lambda_0,\sigma,s}$ is finite for some $s \in (0,1)$, $\lambda_0 > 0$, $\sigma > \frac{d}{2}+\frac{s}{2}+2$ and $M>\frac{d}{2}+1$. Then, the unique classical solution $(f,u) \in C(0,T_{max};H^\sigma_{x,v;M})\times C(0,T_{max};H^\sigma_{x})$ satisfies for all $t \in [0,T_{max})$ the upper bounds
\begin{equation}\label{upper bound1}
\|f\|_{\lambda,\sigma,M,s} \leqslant C_1(1+t)g(t)
\end{equation}
and 
\begin{equation}\label{upper bound2}
\| u \|_{\lambda,\sigma,s} \leqslant \left(\| u_0 \|_{\lambda_0,\sigma,s}+C_2\int_0^t (1+\tau)g(\tau) \ud \tau\right)\exp\bigg[C_2 \int_0^t g(\tau) \ud \tau \bigg] ,
\end{equation}
and for all $t \in [0,T_{max})$ the lower bound
\begin{equation}\label{lower bound} 
\lambda(t) \geqslant  \left(2C_3t+\lambda_0^{-1} \right)^{-1} \exp \bigg[-C_3 \int_0^t \left( 1+ \|u(\tau)\|_\sigma+\| f(\tau) \|_{\sigma,M} \right) \ud \tau \bigg] > 0 ,
\end{equation}
where $$ g(t):=  \exp\big[C_0 \int_0^t \left(\|u(\tau)\|_{W^{1,\infty}}+\| \rho(\tau)\|_\infty+\| f (\tau) \|_{\infty,M}^2+1 \right)\ud \tau\big] , $$
and where $T_{max}$ is the maximal time of existence. The constants $C_0$, $C_1$, $C_2$ and $C_3$ depend on the initial data $(f_0,u_0)$, the radius of regularity $\lambda_0$, the Sobolev correction $\sigma$, the weight $M$ and the dimension $d$.
\end{theorem}

\begin{remark}
We could remove the term $\|\rho(\tau)\|_\infty$ which appears in the definition of $g$ as well as the compact support assumption on the initial data, however the estimates for the growth of the Gevrey norm and the radius of regularity would be bounded in terms of $\|f\|_{\sigma,M}$ and $\|u\|_{\sigma}$, instead of $\|u\|_{W^{1,\infty}}+ \|\rho_f\|_\infty + \| f (t) \|_{\infty,M}^2 + 1$,  and for a \emph{short time},  due to the proof given in Subsection \ref{subsection main results} (proof of Theorem \ref{propagation Analyticite}).  Thus,  the propagation of Gevrey regularity and analyticity, for global solutions, is not guaranteed in this case. 

 If we want to have a propagation for global solutions (for small data), we must have more finite moments (in velocities) and show a propagation of the moments in this case. This is analogous to the results of Pfaffelmoser \cite{pfaffelmoser1992global} and Lions-Perthame \cite{lions1991propagation} for the Vlasov-Poisson system. This last constraint comes from the estimate of the commutator that one needs to control the force term which comes from the Vlasov equation.  In particular, in both cases,  we obtain propagation of analyticity for the Vlasov-Navier-Stokes system in $\mathbb{T}^d\times\RR^d$ and $\RR^d\times\RR^d$.
\end{remark}

\begin{theorem}[Propagation of analyticity]\label{propagation Analyticite}  Let $(f_0,u_0)$ be initial data for the VNS system \eqref{VNS} on $\mathbb{T}^d\times\RR^d$ (or $\RR^d\times\RR^d$) such that $\|f_0\|_{\lambda_0,\sigma,M,1} + \|u_0\|_{\lambda_0,\sigma,1}$ is finite for some $\lambda_0 > 0$, $\sigma > \frac{d}{2}+\frac{5}{2}$ and $M>\frac{d}{2}+1$. Then, the classical solution $(f,u) \in C(0,T_{max};H^\sigma_{x,v;M})\times C(0,T_{max};H^\sigma_{x})$ satisfies for all $t \in [0,T_{max})$ the upper bounds 
\begin{equation}\label{upper bound analytic}
 \|f(t)\|_{\lambda,\sigma,M,1} \leqslant \|f_0\|_{\lambda_0,\sigma,M,1} \ \exp\big[ C_4\int_0^t (1+\|u(\tau)\|_\sigma) \ud \tau \big]
\end{equation}
and 
\begin{equation}\label{upper bound analytic}
\| u \|_{\lambda,\sigma,1} \leqslant \bigg(\|u_0\|_{\lambda_0,\sigma,1} + C_5 \int_0^t \|f(\tau)\|_{\lambda,\sigma,M,1}  \ud \tau \bigg)\exp\bigg[C_5\int_0^t Y(\tau) \ud \tau\bigg],
\end{equation}
and the lower bound \eqref{lower bound}, where $Y(t) := \|u(t)\|_{\sigma} + \|f(t)\|_{\sigma,M}$ and $T_{max}$ is the maximal time of existence. The constants $C_4$ and $C_5$ depend on the initial data $(f_0,u_0)$, the radius of regularity $\lambda_0$, the Sobolev correction $\sigma$, the weight $M$ and the dimension $d$.
\end{theorem}

\begin{remark} 
\ni \begin{enumerate}
\item We assume that $Y(0) := \|u_0\|_{\sigma} + \|f_0\|_{\sigma,M}$ is finite for $\sigma > \frac{d}{2}+2$ and $M> \frac{d}{2}+1$. Then there exists $T > 0$ that depends only on $Y_0$, such that for all $t \in (0,T)$
we have the estimate 
$$ \|u(t)\|_{\sigma} + \|f(t)\|_{\sigma,M} \leqslant \sqrt{2}\frac{(1+\|u_0\|_{\sigma}^2 + \|f_0\|_{\sigma,M}^2)^{\frac{1}{2}}}{1-Ct(1+\|u_0\|_{\sigma}^2 + \|f_0\|_{\sigma,M}^2)^{\frac{1}{2}}} ,$$
where $C$ is a constant that depends on $\sigma$, $M$ and $d$.  The above inequality is established at the beginning of the proof of Theorem \ref{propagation Analyticite}. 
\item Theorem \ref{propagation Analyticite} remains valid for $s \in (0,1)$, which gives us propagation of the Gevrey regularity without any assumption on the support of $f$, but this being for a time interval which depends on the size (smallness) of the initial data. 
\end{enumerate}
\end{remark}

\begin{theorem}[Blow up criterion]\label{Blow up}
Let $(f_0,u_0)$ be initial data for the VNS system \eqref{VNS} on $\mathbb{T}^d\times\RR^d$ (or $\RR^d\times\RR^d$) such that $\|f_0\|_{\lambda_0,\sigma,M,s} + \|u_0\|_{\lambda_0,\sigma,s}$ is finite for some $s \in (0,1)$, $\lambda_0 > 0$,  $\sigma > \frac{d}{2}+\frac{s}{2}+2$ and $M>\frac{d}{2}+1$. Let $T_{\mathrm{max}}$ be the maximal time of existence of the Gevrey solutions $(f,u)$ of the VNS system \eqref{VNS}. Then, if for some $T \in [0,T_{\mathrm{max}}]$,  we have
\begin{equation}
 \int_0^T \left(\|u(\tau)\|_{W^{1,\infty}}+\| \rho(\tau)\|_\infty+\| f(\tau) \|_{\infty,M}^2+1 \right)\ud \tau < +\infty ,  
\end{equation}
then $T < T_{\mathrm{max}}$.  Similarly, for the case $s=1$, if we have for some $T \in [0,T_{\mathrm{max}}]$
\begin{equation}
 \int_0^T \left( \|f(\tau)\|_{\sigma,M} + \|u(\tau)\|_{\sigma} \right)\ud \tau < +\infty , 
\end{equation}
then $T < T_{\mathrm{max}}$.
\end{theorem}

In other words, the propagation of Gevrey regularity on $[0,T]$ follows as long as $\|f\|_{\sigma,M} + \|u\|_{\sigma}$ is  bounded in $L^1(0,T)$. 

As an application of Theorem \ref{propagation Gevrey}, we obtain global existence of Gevrey solutions for the VNS  system \eqref{VNS} in $\mathbb{T}^3\times\RR^3$ for small data. This result follows directly by using the results of Han-Kwan, Moussa and Moyano in \cite{HMM}.

\begin{corollary}[Global existence of Gevrey solutions for $s \in (0, 1)$]\label{existence globale}
 Let $(f_0,u_0)$ be initial data for the VNS system \eqref{VNS} in $\mathbb{T}^3\times\RR^3$ such that, $f_0$ has  a compact support in velocity and $\|f_0\|_{\lambda_0,\sigma,M,s} + \|u_0\|_{\lambda_0,\sigma,s}$ is finite  for some $s \in (0,1)$, $\lambda_0 > 0$, $\sigma > \frac{7}{2}+\frac{s}{2}$ and $M > 4$.  Let $\mathcal{E}(0)$ (defined in \eqref{energie modulee}) small enough in the sense of Theorem \ref{thm HMM}. Then, there exist a unique global classical solution $(f,u) \in C (0, \infty;H^\sigma(\mathbb{T}^3\times\RR^3))\times C(0, \infty;H^\sigma(\mathbb{T}^3))\cap L^2(0,\infty;H^{\sigma+1}(\mathbb{T}^3))$ of the VNS system \eqref{VNS} satisfies for all $t \in [0,\infty)$ the upper bounds \eqref{upper bound1} and \eqref{upper bound2}, and the radius of regularity $\lambda(t)$ satisfies  the lower bound \eqref{lower bound}. 
\end{corollary}

\ni \textbf{Comments.  }
\begin{enumerate}
\item Note that we have recovered the same Gevrey estimates for the Vlasov solution as those for Vlasov-Poisson \cite{VR}, 
but with a lower bound for $\lambda$ as $\exp(-\exp(\exp t))$, because of the term $\na_v \cdot(v f)$ which induces $\exp (t)$ in the support of $f$ at time $t$.
\item In \cite{FoiTem}, the radius of analyticity for the Navier-Stokes equations is given by $\lambda(t)=\min(t,\lambda_1,T^*)$, where $\lambda_1$ is the radius of analyticity of the source term $F$ (and which does not depend on $u$) and $T^*$ is the maximal time of existence which depends on $u_0$ and the source term $F$. This implies that at the time $t=0$, $\lambda(0)=0$.  Then, we recover the Sobolev norm of $u_0$ instead of a Gevrey norm at $t=0$.  Thus, a control of the Gevrey norm of $u$ at time $t$ by the Gevrey norm at the initial time $t=0$, allowed them to get the Gevrey propagation for $u_0$ just Sobolev.
\end{enumerate}

\subsubsection*{Idea of the proof and outline of the paper}
The main result follow by energy estimates based on a Fourier space method motivated by the approach used in \cite{KV1} to study the propagation of analytic regularity for the $2d$-Euler system and \cite{VR} for the Gevrey regularity for the Vlasov-Poisson system. The Gevrey norm will be estimated by the Sobolev norm, so we will start with the Sobolev estimates for the solution $(f,u)$ in Section \ref{section Sobolev}, then move on to Gevrey estimates in Section \ref{section Gevrey}.
The parameters $s$, $\sigma$ and $M$ are fixed, while $\lambda$ is a function in $t$. The weight in which $M$ intervenes is made to control the moments $\rho_f$ and $j_f$ in term of the density distribution $f$ as we said earlier, and the function $\lambda$ will be chosen so that the norms representing a loss of Sobolev regularity (due to time derivatives in the energy method) are absorbed. 

\ni \textbf{\emph{Acknowledgements. }}  The author would like to thank Iv\'an Moyano for suggesting this subject, and Marjolaine Puel for her suggestions during the first reading of the paper. Finally, the author would also like to thank the anonymous referees for their several suggestions to improve the article. 

\section{Sobolev estimates}\label{section Sobolev}
The purpose of this section is to show the following proposition on the propagation of the Sobolev regularity. 
For this purpose, we prove a quantitative bound for the growth of weighted Sobolev norms of the solutions $(u,f)$ of the VNS system in terms of $\| u \|_{W^{1,\infty}}$,  $\|\rho\|_\infty$ and the support of $f$ in velocity $V^M(t)$.  
\begin{proposition}[Sobolev estimates for VNS]\label{Sobolev VNS}
 Let $\sigma>0$ and $M>\frac{d}{2}+1$. Let $(f,u)$ satisfying $\eqref{VNS}$ such that $f$ has a compact support in velocity and $\|f_0\|_{\sigma,M}+\|u_0\|_\sigma $ is finite. Then, the following estimate holds
\begin{equation}\label{Sobolev-VNS2}
\|f(t)\|_{\sigma,M}^2 + \| u(t)\|_\sigma^2  \leqslant (\|f_0\|_{\sigma,M}^2 + \| u_0\|_\sigma^2)g(t)
\end{equation}
where $$ g(t):=  \exp\big[C_0 \int_0^t \left(\|u(\tau)\|_{W^{1,\infty}}+\| \rho(\tau)\|_\infty+\| f(\tau) \|_{\infty,M}^2+1 \right)\ud \tau\big] , $$
and $C_0$ is a positive constant which depends only on $\sigma$,  $M$ and $d$. 
\end{proposition} 
 
 In order to prove Proposition \ref{Sobolev VNS}, we will establish estimates on the time derivative of the Sobolev norm of each of the Vlasov and Navier-Stokes solutions, then we conclude with Gronwall's Lemma applied to a combination of the two estimates. 
 
\subsection{Sobolev estimates for solutions of Vlasov's equation}
The goal of this subsection is to prove the following lemma:
\begin{lemma}[Sobolev estimates for Vlasov]\label{lem Sob V}
Let $\sigma>0$ and $M>0$. Let $(u,f)$ be the solution of $\eqref{VNS}$ such that $f$ has a compact support in velocity.  Then, one has the following estimate 
\begin{equation}\label{Sobolev-V2}
\frac{\ud }{\ud t}\| f \|_{\sigma,M}^2\lesssim (\| u\|_{W^{1,\infty}}+1)\| f\|_{\sigma,M}^2+ \| f \|_{\infty,M} \| u\|_{\sigma+1} \| f \|_{\sigma,M}.
\end{equation}
\end{lemma}

\ni We will start with two lemmas. The first one is on the Gevrey norm estimate for the density $\rho_f$ and the moment $j_f$, which gives in particular the Sobolev estimates for $\lambda=0$. 

\begin{lemma}[Density and first moment estimates]\label{rho,j < f}
 Let $s\in(0,1], \lambda \geqslant 0$ and $\sigma>0$. Let $f\in \mathcal{G}^{\lambda,\sigma,M,\frac{1}{s}}(\mathbb{T}^d\times\RR^d)$. Then, for $M > \frac{d}{2} + 1$, one has  
\begin{equation}\label{rho < f}
\| \rho_f \|_{\lambda,\sigma,s}  \lesssim \| f \|_{\lambda,\sigma,M,s}  \qquad \mbox{and} \qquad \| j_f \|_{\lambda,\sigma,s}  \lesssim \| f \|_{\lambda,\sigma,M,s}.   
\end{equation}
In particular, for $\lambda = 0$ we get the following Sobolev estimates:
$$ \| \rho_f \|_{\sigma}  \lesssim \| f \|_{\sigma,M}  \qquad \mbox{and} \qquad \| j_f \|_{\sigma}  \lesssim \| f \|_{\sigma,M}.   $$
\end{lemma}
\begin{remark}
The inequality 
$$\| \rho_f \|_{\lambda,\sigma,s}  \lesssim \| f \|_{\lambda,\sigma,M,s} $$
is valid as soon as $M > \frac{d}{2}$.
\end{remark}

\ni \bp The proof of Lemma \ref{rho,j < f} is obtained by the Cauchy-Schwarz inequality.  Another proof concerning the inequality on $\rho_f$ is given in \cite{VR}. \\
We have:
\begin{align*}
\| \rho_f\|_{\lambda,\sigma,s}^2&=\|\Lambda^\sigma \ue^{\lambda\Lambda^s}\rho_f\|_{L^2}^2=\int_{\mathbb{T}^d} \bigg|\int_{\mathbb{R}^d}\Lambda^\sigma \ue^{\lambda\Lambda^s}f\ud v\bigg|^2\ud  x \\
&\leqslant \bigg(\int_{\mathbb{R}^d}\big(1+|v|\big)^{-2M}\ud v\bigg)\iint_{\mathbb{T}^d\times\mathbb{R}^d} \big|\big(1+|v|\big)^M\Lambda^\sigma \ue^{\lambda\Lambda^s}f\big|^2\ud v \ud x \\
&\lesssim \|f\|_{\lambda,\sigma,M,s}^2 \quad \mbox{ for } M> \frac{d}{2} .
\end{align*}
Similarly for $j_f$,  we have:
\begin{align*}
\| j_f\|_{\lambda,\sigma,s}^2&=\int_{\mathbb{T}^d} \bigg|\int_{\mathbb{R}^d}v\Lambda^\sigma \ue^{\lambda\Lambda^s}f\ud v\bigg|^2\ud x \\
&\leqslant \bigg(\int_{\mathbb{R}^d}\big(1+|v|\big)^{2-2M}\ud v\bigg)\iint_{\mathbb{T}^d\times\mathbb{R}^d} \big|\big(1+|v|\big)^{M}\Lambda^\sigma \ue^{\lambda\Lambda^s}f\big|^2\ud v\ud x \\
&\lesssim \|f\|_{\lambda,\sigma,M,s}^2 \quad \mbox{ for } M> \frac{d}{2} +1 .
\end{align*}
\ep

\ni The second lemma concerns the estimation of some products of functions and commutators. 
The first two inequalities of the lemma were taken from \cite[Claim 2]{VR}, whose proof was inspired by \cite{klainerman1981singular}, while the last two inequalities were taken from the last reference.
\begin{lemma}[Inequalities on some products and commutators]\label{commutateur} 
\item Let $u, v\in W^{1,\infty}(\mathbb{T}^d)\bigcap H^\sigma(\mathbb{T}^d)$ and let $f\in L^\infty(\mathbb{T}^d\times\mathbb{R}^d)\bigcap H^\sigma(\mathbb{T}^d\times\mathbb{R}^d)$  be a function with compact support in velocity
$$ \supp f \subset \mathbb{T}^d \times \mathrm{B}(0,C_f) ,  \quad  C_f \in (0,+\infty) .$$
Then, the following inequalities hold
\begin{enumerate}
\item $\ds \| uf \|_\sigma \lesssim_{C_f} \| u \|_\infty\| f\|_\sigma + \| u\|_\sigma \| f\|_\infty.$
\item $\ds \underset{|\beta|\leqslant \sigma}{\sum}  \| \uD^\beta(uf)-u\uD^\beta(f) \|_2 \lesssim_{C_f} \| \nabla u \|_\infty\| f\|_{\sigma-1} + \| u\|_\sigma \| f\|_\infty.$ 
\item $\ds \|uv\|_\sigma \lesssim \|\na_x u\|_\infty \|v\|_{\sigma-1} +\|u\|_{\sigma}  \|v\|_{\infty}  .  $ 
\item $\ds \underset{|\beta|\leqslant \sigma}{\sum}   \| \uD^\beta(uv)-u\uD^\beta(v) \|_2 \lesssim \| \nabla_x u \|_\infty\| v\|_{\sigma-1} + \| u\|_\sigma \| v\|_\infty.$
\end{enumerate}
\end{lemma}

\begin{remark}
The condition of $f$ with compact support eliminates the analytic case, i.e. $s=1$.  Throughout the paper, the assumption on the support of $f$ is only used to obtain inequalities 1 and 2 of the above lemma, which are used to obtain the Sobolev estimates of the Vlasov solution (Lemma \ref{estimation Sobolev pour f}).  This means that the analytic case is treated separately (proof of Theorem \ref{propagation Analyticite}), in order to control the Sobolev norms only. We refer to \cite[Proof of Lemma 5.2.5]{dechicha2023fractional} for a detailed proof of the above lemma, to see where the support assumption comes into play more precisely.
\end{remark}

\begin{remark}
One expects that the condition ``$f$ with compact support'' can be replaced by a weight for the velocity variable in order to be able to apply Jensen's inequality on $\RR^d$ with the Lebesgue measure, which plays a crucial role in the proof,  but this requires more finite moments on $f$ and which amounts to proving a propagation of moments (see e.g. \cite[Lemma 1]{HKR2016}).
\end{remark}

\ni \bpl \ref{lem Sob V}.  The proof of this lemma is based on that of Lemma 1 in \cite{VR}.  Recall that
$$
\|f\|_{\sigma,M}^2:=\sum_{|\alpha|\leqslant M} \|v^\alpha f\|_\sigma^2 = \sum_{|\beta|\leqslant \sigma}\sum_{|\alpha|\leqslant M} \|\uD^\beta v^\alpha f\|_{L^2}^2 .
$$
We will show that for all $|\beta| \leqslant \sigma$ and $|\alpha| \leqslant M$,
\begin{align*}
 \frac{1}{2}\frac{\ud}{\ud t} \|\uD^\beta (v^\alpha f)\|_{L^2}^2 &\lesssim (\|u\|_{W^{1,\infty}}+1)\bigg[ \|\uD^\beta (v^\alpha f)\|_{L^2}^2 + \|\uD^\beta (v^{\alpha-e_i} f)\|_{L^2}^2 \bigg] + \|f\|_{\infty,M}\|u\|_{\sigma+1}\|\uD^\beta (v^\alpha f)\|_{L^2},
\end{align*}
where $e_i$ denote the multi-index worth $1$ in the i-th position and $0$ elsewhere, i.e. $e_i :=(0,\ldots,1,\ldots,0)$.  The estimate \eqref{Sobolev-V2} is obtained after summation over $\alpha$ and $\beta$ in the previous inequality. By using the Vlasov equation, we write
\begin{align*}
\frac{1}{2}\frac{d}{\ud t} \|\uD^\beta (v^\alpha f)\|_{L^2}^2 &= \int_{\mathbb{T}^d\times\mathbb{R}^d} \uD^\beta (v^\alpha f) \uD^\beta (v^\alpha \partial_t f) \ud x\ud v \\
&= - \int_{\mathbb{T}^d\times\mathbb{R}^d} \uD^\beta (v^\alpha f) \uD^\beta (v^\alpha v\cdot \nabla_x f) \ud x\ud v \\
& \ \ \ - \int_{\mathbb{T}^d\times\mathbb{R}^d} \uD^\beta (v^\alpha f) \uD^\beta \big(v^\alpha  \nabla_v \cdot\big[(u-v)f\big]\big) \ud x\ud v \\
&=: - (E+F) ,
\end{align*}
with 
$$ E := \int_{\mathbb{T}^d\times\mathbb{R}^d} \uD^\beta (v^\alpha f) \uD^\beta (v^\alpha v\cdot \nabla_x f) \ud x\ud v \quad \mbox{ and } \quad F := \int_{\mathbb{T}^d\times\mathbb{R}^d} \uD^\beta (v^\alpha f) \uD^\beta \big(v^\alpha  \nabla_v \cdot\big[(u-v)f\big]\big) \ud x\ud v .$$
\textbf{Estimation of $E$.} We have
\begin{align*}
E &=  \sum_{i=1}^d  \int_{\mathbb{T}^d\times\mathbb{R}^d} \uD^\beta (v^\alpha f) \partial_{x_i}\big(\uD^{\beta} [v_i(v^\alpha f)]\big) \ud x\ud v \\
&= \sum_{i=1}^d  \int_{\mathbb{T}^d\times\mathbb{R}^d} \uD^\beta (v^\alpha f) \sum_{\beta_1+\beta_2=\beta} \begin{pmatrix} \beta \\ \beta_1  \end{pmatrix}  \uD^{\beta_1}(v_i)\partial_{x_i} \big(\uD^{\beta_2} (v^\alpha f)\big) \ud x\ud v ,
\end{align*}
where we used Leibniz's formula in the last line,  and where each term in this line is zero for $|\beta_1|>1$. Moreover,  for the case $|\beta_1| \leqslant 1$, either $\uD^{\beta_1}(v_i)=0$, $\uD^{\beta_1}(v_i)=1$ or $\uD^{\beta_1}(v_i)=v_i$.  Then,
$$
E= \sum_{i=1}^d  \int \bigg[ v_i\uD^\beta (v^\alpha f) \partial_{x_i} \big(\uD^{\beta} (v^\alpha f)\big) + \uD^\beta (v^\alpha f) \partial_{x_i} \big(\uD^{\beta-e_i} (v^\alpha f)\big) \bigg]\ud x\ud v := E_1+E_2 ,
$$
where 
$$ E_1:= \sum_{i=1}^d  \int_{\mathbb{T}^d\times\mathbb{R}^d} v_i\uD^\beta (v^\alpha f) \partial_{x_i} \big(\uD^{\beta} (v^\alpha f)\big) \ud x\ud v = \frac{1}{2} \sum_{i=1}^d  \int_{\mathbb{T}^d\times\mathbb{R}^d}  v_i\partial_{x_i}|\uD^\beta (v^\alpha f)|^2\ud x\ud v =0 
$$
and 
\begin{align*}
|E_2|&:= \bigg|\sum_{i=1}^d  \int_{\mathbb{T}^d\times\mathbb{R}^d}  \uD^\beta (v^\alpha f) \partial_{x_i} \big(\uD^{\beta-e_i} (v^\alpha f)\big) \ud x\ud v \bigg| \\
&\leqslant \frac{1}{2} \sum_{i=1}^d \bigg(\int_{\mathbb{T}^d\times\mathbb{R}^d}  |\uD^\beta (v^\alpha f)|^2\ud x\ud v + \int_{\mathbb{T}^d\times\mathbb{R}^d}  \big|\partial_{x_i} \big(\uD^{\beta-e_i} (v^\alpha f)\big)\big|^2 \ud x\ud v\bigg) \lesssim \|\uD^\beta(v^\alpha f)\|_{L^2}^2 .
\end{align*}
\textbf{Estimation of $F$.} We will expand the scalar product and make the commutator appear in order to use inequality 2 of Lemma \ref{commutateur} as follows:
\begin{align*}
F&:=  \int_{\mathbb{T}^d\times\mathbb{R}^d} \uD^\beta (v^\alpha f) \uD^\beta \big(v^\alpha  \nabla_v \cdot\big[(u-v)f\big]\big) \ud x\ud v \\
&= \sum_{i=1}^d \int_{\mathbb{T}^d\times\mathbb{R}^d} \uD^\beta (v^\alpha f) \uD^\beta \big(v^\alpha  \partial_{v_i}\big[(u_i-v_i)f\big]\big) \ud x\ud v \\
&=\sum_{i=1}^d \int_{\mathbb{T}^d\times\mathbb{R}^d} \uD^\beta (v^\alpha f)\bigg[ \uD^\beta (u_i v^\alpha  \partial_{v_i} f)-\uD^\beta(v^\alpha f)-\uD^\beta(v^\alpha v_i\partial_{v_i}f)\bigg] \ud x\ud v \\
&=\sum_{i=1}^d \int_{\mathbb{T}^d\times\mathbb{R}^d} \uD^\beta (v^\alpha f)\bigg[ \partial_{v_i}\uD^\beta (u_i v^\alpha f)-u_i \partial_{v_i}\uD^\beta(v^\alpha f) + u_i \partial_{v_i}\uD^\beta(v^{\alpha} f) -\uD^\beta(u_i v^{\alpha-e_i} f)   \\ 
& \hspace{3.5cm} -\uD^\beta(v^\alpha f) -\uD^\beta(v^\alpha v_i\partial_{v_i}f)\bigg] \ud x\ud v \\
&= F_1 +F_2- F_3 - d \|\uD^\beta(v^\alpha f)\|_{L^2}^2 - F_4 ,
\end{align*}
where 
$$ F_1:=\sum_{i=1}^d \int \uD^\beta (v^\alpha f)\big[  \partial_{v_i} \uD^\beta (u_i v^\alpha f)-u_i  \partial_{v_i} \uD^\beta (v^\alpha  f)\big]\ud x\ud v , \quad  F_2:= \sum_{i=1}^d \int  u_i \uD^\beta (v^\alpha f) \partial_{v_i}\uD^\beta (v^\alpha f) \ud x\ud v , $$
$$  F_3:= \sum_{i=1}^d \int  \uD^\beta (v^\alpha f) \uD^\beta (u_i  v^{\alpha-e_i} f) \ud x\ud v \quad \mbox{ and } \quad F_4:= \sum_{i=1}^d \int \uD^\beta (v^\alpha f) \uD^\beta(v^\alpha v_i\partial_{v_i}f) \ud x\ud v .$$
For $F_1$, by inequality 2 of Lemma \ref{commutateur}, we obtain
\begin{align*}
 |F_1| &\lesssim \sum_{i=1}^d \| \uD^\beta (v^\alpha f)\|_{L^2} \big( \|\nabla u_i\|_\infty \|v^\alpha f\|_{\sigma}+\|u_i\|_{\sigma+1} \|v^\alpha f\|_\infty \big)  \\
&\lesssim  \| \uD^\beta (v^\alpha f)\|_{L^2} \big( \|\nabla u\|_\infty \|f\|_{\sigma,M}+\|u\|_{\sigma+1} \|f\|_{\infty,M}\big) .
\end{align*}
For $F_2$, we have
$$ F_2 = \frac{1}{2}\sum_{i=1}^d \int  u_i \partial_{v_i}|\uD^\beta (v^\alpha f)|^2 \ud x\ud v = 0 .$$
The term $F_3$, dealt with the same way as $F_1$, so we write
\begin{align*}
 F_3 &= \sum_{i=1}^d \int \uD^\beta (v^\alpha f) \big[ \uD^\beta(u_i v^{\alpha-e_i}f)- u_i \uD^\beta (v^{\alpha-e_i} f)+ u_i \uD^\beta (v^{\alpha-e_i} f)\big]\ud x\ud v  \\
&\lesssim  \| \uD^\beta (v^\alpha f)\|_{L^2} \big( \|\nabla u_i\|_\infty \|v^\alpha f\|_{\sigma-1}+\|u_i\|_{\sigma} \|v^\alpha f\|_\infty + \|u_i\|_\infty \|\uD^\beta (v^{\alpha-e_i} f)\|_{L^2}\big)  \\
&\lesssim  \| \uD^\beta (v^\alpha f)\|_{L^2} \big( \|\nabla u\|_\infty \|f\|_{\sigma,M}+\|u\|_{\sigma} \|f\|_{\infty,M}+\|u\|_\infty \|\uD^\beta (v^{\alpha-e_i} f)\|_{L^2}\big) .
\end{align*}
Finally, for $F_4$, we have
\begin{align*}
F_4 &= \sum_{i=1}^d \sum_{\beta_1+\beta_2=\beta} \begin{pmatrix} \beta \\ \beta_1  \end{pmatrix}  \int \uD^\beta (v^\alpha f) \uD^{\beta_1}(v_i)\partial_{v_i}\uD^{\beta_2}(v^\alpha f) \ud x\ud v \\
&= \sum_{i=1}^d \int \uD^\beta (v^\alpha f) \big[v_i \partial_{v_i}\uD^{\beta}(v^\alpha f)+\partial_{v_i}\uD^{\beta-e_i}(v^\alpha f) \big] \ud x\ud v \\
&= \sum_{i=1}^d \int  \frac{v_i }{2}\partial_{v_i}|\uD^{\beta}(v^\alpha f)|^2+\sum_{i=1}^d \int \uD^\beta (v^\alpha f) \partial_{v_i}\uD^{\beta-e_i}(v^\alpha f) \big] \ud x\ud v  ,
\end{align*}
which implies that 
$$
|F_4| \leqslant \frac{1}{2} \sum_{i=1}^d  \int \bigg[ 2|\uD^{\beta}(v^\alpha f)|^2+ |\partial_{v_i}\uD^{\beta-e_i}(v^\alpha f)|^2 \bigg] \ud x\ud v \lesssim \|\uD^\beta (v^\alpha f)\|_{L^2}^2 .
$$
Thus, by combining the inequalities on $E_i$ and $F_i$, we obtain
\begin{align*}
 \frac{1}{2}\frac{\ud}{\ud t} \|\uD^\beta (v^\alpha f)\|_{L^2}^2 \lesssim &(\|u\|_{W^{1,\infty}}+1)\bigg[ \|\uD^\beta (v^\alpha f)\|_{L^2}^2 + \|\uD^\beta (v^{\alpha-e_i} f)\|_{L^2}^2 \bigg] \\
 &+ \|f\|_{\infty,M}\|u\|_{\sigma+1}\|\uD^\beta (v^\alpha f)\|_{L^2},
\end{align*}
and by summing over $\alpha$ and $\beta$, we obtain for $\sigma>0$ and $M>0$
$$
 \frac{1}{2}\frac{\ud}{\ud t} \|f\|_{\sigma,M}^2 \lesssim (\|u\|_{W^{1,\infty}}+1)\|f\|_{\sigma,M}^2 + \|f\|_{\infty,M}\|u\|_{\sigma+1}\|f\|_{\sigma,M}  .
$$
\epl

\subsection{Sobolev estimates for  Navier-Stokes}
In this subsection, we estimate the Sobolev norm of the solution of the Navier-Stokes equations, which we summarize in the following lemma:
\begin{lemma}[Sobolev estimates for NS]
Let $\sigma>0$ and $M>\frac{d}{2}+1$. Let $(f,u)$ satisfying equations \eqref{VNS}. Then, one has the following estimate
\begin{equation}\label{Sobolev-NS2}
\frac{1}{2}\frac{\ud}{\ud t}\| u\|_\sigma^2 + \| u\|_{\sigma+1}^2 \lesssim (\| \na_x u\|_\infty+\|\rho\|_\infty)\| u\|_\sigma^2 + (\| u\|_\infty+1)\| u\|_\sigma \| f\|_{\sigma,M} .
\end{equation}
\end{lemma}

\ni \bp Recall that $u$ satisfies the equations 
$$ \pa_t u + u \cdot \na_x u - \Delta_x u + \na_x p = j_f - \rho_f u \quad \mbox{ and } \quad \na_x \cdot u = 0 ,$$
and one has:
$$ \frac{1}{2} \frac{\ud}{\ud t} \| u \|_\sigma^2 = \frac{1}{2} \sum_{|\alpha|\leqslant \sigma} \frac{\ud}{\ud t} \| \uD^\alpha u\|_{L^2}^2 .$$
Then, applying $\uD^\alpha$ to the first equation of $u$ and integrating it against $\uD^\alpha u$, we get:
$$ \frac{\ud}{2\ud t} \| \uD^\alpha u\|_{L^2}^2 - \langle \Delta_x(\uD^\alpha u),\uD^\alpha u \rangle_{L^2} = -\langle \uD^\alpha (u\cdot \na_x u),\uD^\alpha u \rangle_{L^2} + \langle \uD^\alpha j,\uD^\alpha u \rangle_{L^2} - \langle \uD^\alpha (\rho u),\uD^\alpha u \rangle_{L^2} .$$
Since $\na_x \cdot u = 0$ then,  $ \langle u\cdot \na_x(\uD^\alpha u),\uD^\alpha u \rangle_{L^2} = 0$. Therefore, 
\begin{align*}
\big| \langle \uD^\alpha (u\cdot \na_x u),\uD^\alpha u \rangle_{L^2} \big| &= \big| \langle \uD^\alpha (u\cdot \na_x u) - u\cdot \na_x(\uD^\alpha u),\uD^\alpha u \rangle_{L^2} \big|  \\
&\leqslant \| \uD^\alpha u \|_{L^2}  \big\| \uD^\alpha (u\cdot \na_x u) - u\cdot \na_x(\uD^\alpha u) \big\|_{L^2}.
\end{align*}
Thus, by Lemma \ref{commutateur} we obtain:
$$ \big| \langle \uD^\alpha (u\cdot \na_x u),\uD^\alpha u \rangle_{L^2} \big|  \lesssim \| \na_x u \|_\infty \| \uD^\alpha u\|_{L^2}^2 $$
and 
$$ \| \uD^\alpha (\rho u)\|_{L^2} \lesssim  \|\rho \|_\infty \|\uD^\alpha u\|_{L^2} + \| \uD^\alpha \rho \|_{L^2} \| u\|_\infty .$$
Therefore,
$$ \frac{1}{2} \frac{\ud}{\ud t} \| u \|_\sigma^2 + \| u \|_{\sigma+1}^2 \lesssim \big(\| \na_x u \|_\infty + \|\rho \|_\infty \big) \|u\|_\sigma^2 + \big(\|u\|_\infty \|\rho\|_\sigma+\|j \|_\sigma\big) \|u\|_\sigma . $$
Finally, thanks to inequalities \eqref{rho < f} of Lemma \ref{rho,j < f} (for $\lambda=0$), we get for $M>\frac{d}{2} +1$:
$$ \frac{1}{2} \frac{\ud}{\ud t} \| u \|_\sigma^2 + \| u \|_{\sigma+1}^2 \lesssim \big(\| \na_x u \|_\infty + \|\rho \|_\infty \big) \|u\|_\sigma^2 + \big(\|u\|_\infty+1\big) \|f\|_{\sigma,M}\|u\|_\sigma . $$
Hence the inequality \eqref{Sobolev-NS2} holds true.
\ep
\subsection{Sobolev estimates for VNS and proof of Proposition \ref{Sobolev VNS}}
In this subsection, we will combine the estimates obtained on the solutions of the Vlasov and Navier-Stokes equations, as we said above, in order to be able to apply Gronwal's lemma and control the Sobolev norm of $(f,u)$ at time $t$ by that at initial time $t=0$.  
\begin{lemma}[Sobolev estimates for VNS]\label{lemme Gronwall Sob}
Let $(f,u)$ satisfying equations \eqref{VNS} such that $f$ has a compact support in velocity. Let $\sigma>0$ and $M>\frac{d}{2}+1$. Then, one has the following estimate
\begin{equation}\label{Gronwall Sobolev}
\frac{\ud}{\ud t} \left( \| u \|_\sigma^2 + \|f\|_{\sigma,M}^2\right) \lesssim \left(\|u\|_{W^{1,\infty}}+\|\rho\|_\infty+\| f\|_{\infty,M}^2+1\right) \left(\|u\|_\sigma^2 + \|f\|_{\sigma,M}^2\right).
\end{equation} 
\end{lemma} 

\ni \bpl \ref{lemme Gronwall Sob}. We obtain \eqref{Gronwall Sobolev} by combining the two inequalities \eqref{Sobolev-V2} and \eqref{Sobolev-NS2} and using Young's inequality.  Indeed,  summing \eqref{Sobolev-V2} and \eqref{Sobolev-NS2} gives:
\begin{align*}
\frac{1}{2} \frac{\ud}{\ud t} \left( \| u \|_\sigma^2 + \|f\|_{\sigma,M}^2\right) +\|u\|_{\sigma+1}^2  & \lesssim \left(\|u\|_{W^{1,\infty}}+1\right) \|f\|_{\sigma,M}^2 + \| f\|_{\infty,M}\|u\|_{\sigma+1}\|f\|_{\sigma,M}  \\
& \ + \left( \| \na_x u \|_\infty +  \|\rho\|_\infty \right) \| u \|_\sigma^2 + \left( \| u \|_\infty + 1 \right) \| u \|_\sigma  \|f\|_{\sigma,M} . 
\end{align*}
Now we need the dissipative term of the NS equations to absorb the term $\|u\|_{\sigma+1}$ that appears on the right-hand side of the inequality.  Thus, by Young's inequality, used for the two terms $\| f \|_{\infty,M}\|u\|_{\sigma+1}\|f\|_{\sigma,M}$ and $\| u \|_\sigma \|f\|_{\sigma,M}$, we obtain: 
\begin{align*}
\frac{1}{2} \frac{\ud}{\ud t} \left( \| u \|_\sigma^2 + \|f\|_{\sigma,M}^2\right) +\|u\|_{\sigma+1}^2 &\leqslant C \left(\|u\|_{W^{1,\infty}}+\|\rho\|_\infty+\| f \|_{\infty,M}^2+1\right) \left(\|u\|_\sigma^2 + \|f\|_{\sigma,M}^2\right) \\
& \ + \frac{1}{2} \|u\|_{\sigma+1}^2.
\end{align*}
 \epl

\ni \bpp \ref{Sobolev VNS}. Inequality \eqref{Sobolev-VNS2} follows from  \eqref{Gronwall Sobolev} and Gronwall's lemma.
\epp

\section{Gevrey estimates for VNS}\label{section Gevrey}
The goal of this section is to prove Theorem \ref{propagation Gevrey}. For this purpose, we will establish quantitative estimates on the Gevrey norms of the Vlasov and Navier-stokes solutions respectively, using a Fourier space
method in analogy to the results proved for the 2D-Euler system in \cite{KV1} and \cite{LO}, and for any dimension for the Vlasov-Poisson system in \cite{VR}. 

\subsection{Preliminary inequalities}
This subsection is summarized in two lemmas that we use throughout this section.  The first contains two discrete Young inequalities (we find a variant of this lemma in \cite{BMM}).
\begin{lemma}[Young's inequality]\label{inegalite de Young}
\item Let $f, \langle k,\eta\rangle^\sigma g \in L^2(\mathbb{Z}^d\times\mathbb{R}^d), \langle k\rangle^\sigma r \in L^2(\mathbb{Z}^d)$ and let $\nu, \beta, \gamma \in \mathbb{R}$. Then,
\begin{enumerate}
\item For $\sigma >\frac{d}{2}+\nu$, one has:
\begin{equation}\label{Young 1}
\bigg| \underset{k,l}{\sum}\int_{\RR^d}f_k(\eta)\langle l\rangle^\nu r_l \langle k-l,\eta\rangle^{\sigma+\beta}g_{k-l}(\eta)\ud \eta \bigg| \lesssim \| f \|_{L^2_{k,\eta}} \| \langle k\rangle^\sigma r\|_{L^2_k} \| \langle k,\eta\rangle^{\sigma+\beta} g\|_{L^2_{k,\eta}} .
\end{equation}
\item For $\sigma>\frac{d}{2}+\gamma-\beta$, one has:
\begin{equation}\label{Young 2}
\bigg| \underset{k,l}{\sum}\int_{\RR^d}f_k(\eta)\langle l\rangle^\sigma r_l \langle k-l,\eta\rangle^{\gamma}g_{k-l}(\eta)\ud \eta \bigg| \lesssim \| f \|_{L^2_{k,\eta}} \| \langle k\rangle^\sigma r\|_{L^2_k} \| \langle k,\eta\rangle^{\sigma+\beta} g\|_{L^2_{k,\eta}}  .
\end{equation}
\end{enumerate}
The constant in the two inequalities depends only on $\nu, \beta, \gamma, \sigma$ and $d$.
\end{lemma}

\ni \bpl \ref{inegalite de Young}. 1. Let $\nu, \beta \in \RR$ and let $\sigma>\frac{d}{2}+\nu$.  We write:
\begin{align*}
&\bigg|\sum_{k,l \in \ZZ^d}\int_{\RR^d} f_k(\eta) \langle l \rangle^\nu r_l \langle k-l,\eta \rangle^{\sigma+\beta} g_{k-l}(\eta) \ud \eta\bigg| \\
&\leqslant \sum_{k} \bigg(\int_{\RR^d} |f_k(\eta)|^2 \ud \eta \bigg)^{\frac{1}{2}}\sum_{l} \langle l \rangle^{\nu}|r_l|\bigg(\int_{\RR^d} \langle k-l,\eta \rangle^{2(\sigma+\beta)} |g_{k-l}(\eta)|^2\ud \eta\bigg)^{\frac{1}{2}} \\
&\leqslant   \bigg(\sum_{k}\int_{\RR^d} |f_k(\eta)|^2 \ud \eta \bigg)^{\frac{1}{2}} \left(\sum_{k}\bigg[\sum_{l} \langle l \rangle^{\nu}|r_l|\bigg(\int_{\RR^d} \langle k-l,\eta \rangle^{2(\sigma+\beta)} |g_{k-l}(\eta)|^2\ud \eta\bigg)^{\frac{1}{2}} \bigg]^2\right)^{\frac{1}{2}}.
\end{align*} 
Now, by Young's inequality for convolution, we get
\begin{align*}
\bigg|\sum_{k,l \in \ZZ^d}\int_{\RR^d} f_k(\eta) \langle l \rangle^\nu r_l \langle k-l,\eta \rangle^{\sigma+\beta} g_{k-l}(\eta) \ud \eta\bigg| &\lesssim \|f\|_{L^2_{k,\eta}} \sum_{k} \langle k \rangle^{\nu}|r_k| \bigg(\sum_{k} \int_{\RR^d} \langle k,\eta \rangle^{2(\sigma+\beta)} |g_k(\eta)|^2\ud \eta\bigg)^{\frac{1}{2}} \\
&\leqslant \bigg(\sum_{k} \langle k \rangle^{2\nu-2\sigma}\bigg)^{\frac{1}{2}}\|\langle k \rangle^{\sigma} r\|_{L^2_k} \|f\|_{L^2_{k,\eta}} \|\langle k,\eta \rangle^{\sigma+\beta} g\|_{L^2_{k,\eta}} \\
&\lesssim_{\nu,\sigma,d} \|\langle k \rangle^{\sigma} r\|_{L^2_k} \|f\|_{L^2_{k,\eta}}  \|\langle k,\eta \rangle^{\sigma+\beta} g\|_{L^2_{k,\eta}} ,
\end{align*} 
since $\ds \sum_k \langle k \rangle^{2\nu-2\sigma} \lesssim 1\ $ for $\ \sigma > \frac{d}{2} +\nu$. 

\ni 2.  The proof of \eqref{Young 2} is analogous.
\epl

\ni The second lemma deals with triangular inequalities on Japanese brackets. 
\begin{lemma}[Some triangular inequalities]\label{inegalites triangulaires}
Let $k,l, \eta \in \mathbb{R}^d$. The following inequalities hold:
\begin{enumerate}
\item $\langle k+l\rangle^s \leqslant \langle k \rangle^s + \langle l \rangle^s, \quad \forall s \in (0,1].$
\item $\langle k+l\rangle^s \leqslant 2^{s-1}(\langle k \rangle^s + \langle l \rangle^s), \quad \forall s \geqslant 1.$
\item $|\langle k\rangle^s-\langle l\rangle^s| \leqslant \langle k -l\rangle^s , \quad \forall s \in (0,1].$
\item $|\langle k\rangle^s-\langle l\rangle^s| \lesssim_s \frac{\langle k -l\rangle}{\langle k\rangle^{1-s}+\langle l\rangle^{1-s}} , \quad \forall s \in (0,1].$
\item $|\langle k,\eta\rangle^s-\langle k-l,\eta\rangle^s| \leqslant \langle l\rangle^s , \quad \forall s \in (0,1].$
\item $|\langle k,\eta\rangle^s-\langle k-l,\eta\rangle^s| \leqslant 2^{s-1}(\langle l\rangle^s+\langle k-l,\eta\rangle^s) , \quad \forall s \geqslant1.$
\item $|\langle k,\eta\rangle^s-\langle k-l,\eta\rangle^s| \lesssim_s \frac{\langle l\rangle}{\langle k,\eta\rangle^{1-s}+\langle k-l,\eta\rangle^{1-s}} , \quad \forall s \in (0,1].$
\end{enumerate}
\end{lemma}
For the proof of this lemma, see Appendix \ref{appendice}. \\

Whether for Vlasov or for Navier-Stokes, the estimates that we are going to establish are based on the two previous lemmas with different parameters.  

\subsection{Gevrey estimates for solutions of Vlasov's equation}\label{sec Gevrey Vlasov}
The aim of this subsection is to prove the following proposition:
\begin{proposition}\label{Gevrey Vlasov}
Assume that the radius of regularity $t \mapsto \lambda(t)$ depends smoothly on time.  Let $s\in(0,1]$, $M>0$ and $\sigma > \frac{d}{2}+\frac{s}{2}+2$. Let $f \in \mathcal{G}^{\lambda,\sigma+\frac{s}{2},M,\frac{1}{s}}(\mathbb{T}^d\times\RR^d)$ and let $u \in \mathcal{G}^{\lambda,\sigma+\frac{s}{2},\frac{1}{s}}(\mathbb{T}^d)$. Then, the following estimate holds
\begin{align}\label{estimation Gevrey pour f}
\frac{1}{2}\frac{\ud}{\ud t} \|f\|_{\lambda,\sigma,M,s}^2 &\lesssim \left( \|u\|_{W^{1,\infty}}+1\right) \|f\|_{\lambda,\sigma,M,s}^2 + \|u\|_{\sigma} \|f\|_{\sigma,M} \|f\|_{\lambda,\sigma,M,s}  \nonumber \\
&+  \left( \dot{\lambda} + \lambda (1+ \|u\|_{\sigma}) + \lambda^2 \|u\|_{\lambda,\sigma,s}\right) \|f\|_{\lambda,\sigma+\frac{s}{2},M,s}^2  \nonumber \\
&+ \left( \lambda \|f\|_{\sigma,M} + \lambda^2 \|f\|_{\lambda,\sigma,M,s} \right) \|u\|_{\lambda,\sigma+\frac{s}{2},s} \|f\|_{\lambda,\sigma+\frac{s}{2},M,s} .
\end{align}
where $\dot \lambda$ denotes the derivative of $\lambda$ with respect to $t$.
\end{proposition}

Subsequently, we will choose $\lambda$ so that the norm $\|f\|_{\lambda,\sigma+\frac{s}{2},M,s}$ will be absorbed in order to control the Gevrey norm of $f$.  Indeed, for $\lambda$ such that
\begin{equation}\label{contrainte 1 lambda}
\dot{\lambda} + \lambda (1+ \|f\|_{\sigma,M} +\|u\|_{\sigma}) + \lambda^2 (\|f\|_{\lambda,\sigma,M,s}+\|u\|_{\lambda,\sigma,s}) \leqslant 0 ,
\end{equation}
we get:
$$ \frac{1}{2}\frac{\ud}{\ud t} \|f\|_{\lambda,\sigma,M,s}^2 \lesssim \left( \|u\|_{W^{1,\infty}}+1\right) \|f\|_{\lambda,\sigma,M,s}^2 + \|u\|_{\sigma} \|f\|_{\sigma,M} \|f\|_{\lambda,\sigma,M,s} .$$
Thus, if the norms $\|u\|_{W^{1,\infty}}$ and $\|u\|_{\sigma}\|f\|_{\sigma,M}$ are in $L^1(0,T)$ then, Gronwall's lemma allows us to conclude. 

\begin{remark}
Note that the estimation of the previous proposition does not require the condition that $f$ has compact support in velocity. 
\end{remark}
As a consequence of Proposition \ref{Gevrey Vlasov}, we have the following estimate which is useful in the case $s=1$.
\begin{corollary}\label{corollaire estimation Sobolev}
Let $M>0$ and $\sigma > \frac{d}{2}+2$. Let $f \in H^\sigma_{M}(\mathbb{T}^d\times\RR^d)$ and let $u \in H^{\sigma}(\mathbb{T}^d)$. Then, the following estimate holds
\begin{align}\label{estimation Sobolev pour f}
\frac{1}{2}\frac{\ud}{\ud t} \|f\|_{\sigma,M}^2 &\lesssim \left( \|u\|_{W^{1,\infty}}+1\right) \|f\|_{\sigma,M}^2 + \|u\|_{\sigma} \|f\|_{\sigma,M}^2.
\end{align}
\end{corollary}
\ni \bpc \ref{corollaire estimation Sobolev}. It suffices to take $\lambda=0$ in Proposition \ref{Gevrey Vlasov} as we only look for Sobolev norms in this case.  Note that the constant behind $\lesssim$ in \eqref{estimation Gevrey pour f} does not depend on $\lambda(0)$.
\epc

\ni \bpp \ref{Gevrey Vlasov}.  We will work in Fourier variables in order to simplify the calculations and expressions and to use the inequalities of the Lemma  \ref{inegalites triangulaires}. The Vlasov equation, in Fourier variables, is given by: 
$$ \pa_t \hat f_k(\eta) - k\cdot \na_\eta \hat f_k(\eta) + \eta \cdot \na_\eta \hat f_k(\eta) + \ui \sum_{l \in \ZZ^d} \hat u_l \cdot \eta \hat f_{k-l} (\eta) = 0 .$$
Recall that
$$ \|f\|_{\lambda,\sigma,M,s}^2 := \sum_{|\alpha|\leqslant M} \|A^\sigma v^\alpha f\|_{L^2}^2 = \sum_{|\alpha|\leqslant M} \|A_k^\sigma(\eta) \uD^\alpha_\eta \hat f\|_{L^2}^2 ,$$
where $$ A^\sigma := A_k^\sigma (\eta) := \langle k,\eta \rangle^{\sigma} \ue^{\lambda \langle k,\eta \rangle^{s}} $$
and $\lambda := \lambda(t)$ is a positive function,  by assumption.  First, we have for a complex-valued function $g$:  
$$ \frac{1}{2} \frac{\ud}{\ud t} \| g \|^2_{L^2_\xi} = \frac{1}{2} \frac{\ud}{\ud t} \int |g|^2 \ \ud \xi = \frac{1}{2} \bigg[\int \bar g \pa_t g \ \ud \xi + \int g \pa_t \bar g \ \ud \xi \bigg] = \Re \int \bar g \pa_t g \ \ud \xi .$$
Then,
\begin{align*}
\frac{1}{2}\frac{\ud}{\ud t}\| f \|_{\lambda,\sigma,M,s}^2 = \frac{1}{2} &\sum_{|\alpha|\leqslant M} \sum_{k \in \ZZ^d} \int_{\RR^d} \pa_t\big(\ue^{2\lambda \langle k,\eta \rangle^{s}}\big)\langle k,\eta \rangle^{2\sigma} |\uD^\alpha_\eta \hat f_k(\eta)|^2 \ud \eta \\
+&\Re \bigg[\sum_{|\alpha|\leqslant M} \sum_{k \in \ZZ^d} \int_{\RR^d} \langle k,\eta \rangle^{2\sigma} \overline{\uD^\alpha_\eta \hat f_k(\eta)} \ \uD^\alpha_\eta\big(\pa_t \hat f(\eta)\big)\  \ud \eta \bigg] .
\end{align*}
Which implies, replacing $\pa_t \hat f_k(\eta)$ by $[ k\cdot \na_\eta \hat f_k(\eta) - \eta \cdot \na_\eta \hat f_k(\eta) - \ui \underset{l }{\sum} \hat u_k \cdot \eta \hat f_{k-l} (\eta)]$ in the last term, that:
\begin{equation}\label{d/dt |f|^2 = lambda'|f|^2+E+F+G}
\frac{1}{2}\frac{\ud}{\ud t}\| f \|_{\lambda,\sigma,M,s}^2 = \dot{\lambda} \| f \|_{\lambda,\sigma+\frac{s}{2},M,s}^2 + \hat E -\hat F - \hat G,
\end{equation}
where
\begin{align*}
\hat E &= \Re  \sum_{|\alpha|\leqslant M} \sum_{k \in \ZZ^d} \int_{\RR^d} \langle k,\eta \rangle^{2\sigma} \ue^{2\lambda \langle k,\eta \rangle^{s}} \overline{\uD^\alpha_\eta \hat f_k(\eta)} \ \uD^\alpha_\eta\big(k\cdot\na_\eta \hat f_k(\eta)\big)\  \ud \eta  , \\
\hat F &= \Re \sum_{|\alpha|\leqslant M} \sum_{k \in \ZZ^d} \int_{\RR^d} \langle k,\eta \rangle^{2\sigma} \ue^{2\lambda \langle k,\eta \rangle^{s}} \overline{\uD^\alpha_\eta \hat f_k(\eta)} \ \uD^\alpha_\eta\big(\eta\cdot\na_\eta \hat f_k(\eta)\big)\ \ud \eta , \\
\hat G &= \Re \sum_{|\alpha|\leqslant M} \sum_{k,l \in \ZZ^d} \int_{\RR^d} \ui \langle k,\eta \rangle^{2\sigma} \ue^{2\lambda \langle k,\eta \rangle^{s}} \overline{\uD^\alpha_\eta \hat f_k(\eta)} \  \hat u_l \cdot \uD^\alpha_\eta\big(\eta \hat f_{k-l}(\eta)\big)\  \ud \eta .
\end{align*}
The terms $\hat E$ and $\hat F$ correspond to the linear part of the equation, $\hat E$ comes from the advection term and $\hat F$ from the force field. Both terms are treated by integration by parts. However, the term $\hat G$, corresponding to the nonlinear part, is the most difficult to handle because of the coupling with the velocity field $u$. In order to reduce the power of the Gevrey norms which comes from the product $uf$ with $f$, we make the commutator between $u$ and $f$ appears by using the fact that $u$ does not depend on $v$. We will do the same to treat the advection term for NS, using the fact that $u$ is divergence-free. \\

\ni \textbf{Estimations of $\hat E$.} By expanding $k \cdot \na_\eta$ and integrating by parts with respect to $\eta$, we obtain:
\begin{align*}
\hat E &= \Re \bigg[\sum_{i=1}^d \sum_{\alpha,k}\int_{\RR^d} A_k^\sigma(\eta) \overline{\uD^\alpha_\eta \hat f_k(\eta)} \ A_k^\sigma(\eta) k_i \pa_{\eta_i} \big(\uD^\alpha_\eta \hat f_k(\eta)\big)\  \ud \eta \bigg] \\
&= -\frac{1}{2}\sum_{i=1}^d \sum_{\alpha,k}\int_{\RR^d} k_i \pa_{\eta_i} \big( A_k^\sigma(\eta)^2\big) | \uD^\alpha_\eta \hat f_k(\eta)|^2  \ud \eta .
\end{align*}
Now, since
$$ \pa_{\eta_i} \left( A_k^\sigma(\eta)\right) = \pa_{\eta_i} \left(\langle k,\eta \rangle^{\sigma} \ue^{\lambda \langle k,\eta \rangle^{s}} \right) = \big[\sigma \eta_i \langle k,\eta \rangle^{\sigma-2} + \lambda s \eta_i \langle k,\eta \rangle^{\sigma+s-2}\big]\ue^{\lambda \langle k,\eta \rangle^{s}} .$$
Then,
$$ \hat E = -\sum_{i,\alpha,k}\int_{\RR^d} |\uD^\alpha_\eta \hat f_k(\eta)|^2 \big[\sigma k_i\eta_i \langle k,\eta \rangle^{2\sigma-2} + \lambda s k_i\eta_i \langle k,\eta \rangle^{2\sigma+s-2}\big]\ue^{2\lambda \langle k,\eta \rangle^{s}} \ud \eta . $$
That leads to
\begin{equation}\label{hat E}
\big| \hat E \big| \lesssim \| f \|_{\lambda,\sigma,M,s}^2 + \lambda \| f \|_{\lambda,\sigma+\frac{s}{2},M,s}^2
\end{equation}
\textbf{Estimations of $\hat F$.} We will proceed as in $\hat E$.
\begin{align*}
\hat F &= \Re  \sum_{i,\alpha,k} \int_{\RR^d} [A_k^\sigma(\eta)]^2 \  \overline{\uD^\alpha_\eta \hat f_k(\eta)} \ \uD^\alpha_\eta \big( \eta_i \pa_{\eta_i} \hat f_k(\eta)\big)\  \ud \eta \\
&=  \Re \sum_{i,\alpha,k} \sum_{\alpha_1+\alpha_2=\alpha} \begin{pmatrix} \alpha_1 \\ \alpha_2  \end{pmatrix} \int_{\RR^d} [A_k^\sigma(\eta)]^2\  \overline{\uD^\alpha_\eta \hat f_k(\eta)} \ \uD^{\alpha_1}_\eta(\eta_i)\ \uD^{\alpha_2}_\eta \big( \pa_{\eta_i} \hat f_k(\eta)\big)\  \ud \eta  \\
&= \sum_{i,\alpha,k} \int_{\RR^d} [A_k^\sigma(\eta)]^2\bigg[ \frac{\eta_i}{2} \pa_{\eta_i}|\uD^\alpha_\eta \hat f_k(\eta)|^2 + \Re \bigg( \overline{\uD^\alpha_\eta \hat f_k(\eta)} \ \pa_{\eta_i}\big(\uD^{\alpha-e_i}_\eta \hat f_k(\eta)\big) \bigg) \bigg] \ud \eta  \\
&=: \hat F_1 + \hat F_2.
\end{align*}
For $\hat F_1$, by integrating by parts with respect to $\eta$ we get:
\begin{align}\label{hat F_1}
\big| \hat F_1 \big| &:= \bigg| \sum_{i,\alpha,k} \int_{\RR^d} [A_k^\sigma(\eta)]^2 \frac{\eta_i}{2} \pa_{\eta_i}|\uD^\alpha_\eta \hat f_k(\eta)|^2 \ud \eta \bigg| \nonumber \\
& =  \bigg| \sum_{i,\alpha,k} \int_{\RR^d} \bigg[ -\frac{1}{2}[A_k^\sigma(\eta)]^2 - \eta_i \pa_{\eta_i}\big(A_k^\sigma(\eta)\big) A_k^\sigma(\eta) \bigg] |\uD^\alpha_\eta \hat f_k(\eta)|^2 \ud \eta \bigg| \nonumber \\
&\lesssim \| f \|_{\lambda,\sigma,M,s}^2 + \lambda \| f \|_{\lambda,\sigma+\frac{s}{2},M,s}^2 ,
\end{align} 
by the same token as for $\hat E$. \\

\ni For $\hat F_2$, we write:
\begin{align}\label{hat F_2}
\big| \hat F_2 \big| &:= \bigg| \Re \sum_{i,\alpha,k} \int_{\RR^d} [A_k^\sigma(\eta)]^2 \ \overline{\uD^\alpha_\eta \hat f_k(\eta)} \ \pa_{\eta_i}\big(\uD^{\alpha-e_i}_\eta \hat f_k(\eta)\big) \ud \eta \bigg| \nonumber \\
& \leqslant \frac{1}{2} \sum_{i,\alpha,k} \int_{\RR^d} [A_k^\sigma(\eta)]^2 \bigg[ |\uD^\alpha_\eta \hat f_k(\eta)|^2 + \big|\pa_{\eta_i}\big(\uD^{\alpha-e_i}_\eta \hat f_k(\eta)\big|^2 \bigg] \ud \eta \nonumber \\
&\lesssim \| f \|_{\lambda,\sigma,M,s}^2  .
\end{align}
Thus, from \eqref{hat F_1} and \eqref{hat F_2} we get:
\begin{equation}\label{hat F}
|\hat F| \lesssim \| f \|_{\lambda,\sigma,M,s}^2+ \lambda \| f \|_{\lambda,\sigma+\frac{s}{2},M,s}^2 .
\end{equation}
\begin{remark}
Note that the only way to absorb the term $\| f \|_{\lambda,\sigma+\frac{s}{2},M,s}$,  representing a possible loss of Sobolev regularity, is to choose a suitable function $\lambda$. We already have $(\dot{\lambda} + \lambda) \| f \|_{\lambda,\sigma+\frac{s}{2},M,s}^2$ which comes from the first three terms of $\frac{\ud}{\ud t}\| f \|_{\lambda, \sigma+\frac{s}{2},M,s}^2$ and others will come from the term $\hat G$.  For this reason we require the loss in the radius of regularity $\lambda(t)$ at an exponential rate. 
\end{remark}

\ni \textbf{Estimations of $\hat G$.} As in $\hat F$, expanding $\hat u_l \cdot \eta$ and using Leibniz, we write:
\begin{align*}
\hat G &= \Re \sum_{j=1}^d \sum_{|\alpha|\leqslant M} \sum_{k,l \in \ZZ^d} \int_{\RR^d} \ui [A_k^\sigma(\eta)]^2\ \overline{\uD^\alpha_\eta \hat f_k(\eta)} \ \uD^\alpha_\eta\big(\hat u_l^j \eta_j \hat f_{k-l}(\eta)\big)\  \ud \eta \\
&= \Imm \sum_{j,\alpha,k,l} \sum_{\alpha_1+\alpha_2=\alpha} \begin{pmatrix} \alpha_1 \\ \alpha_2  \end{pmatrix} \int_{\RR^d} [A_k^\sigma(\eta)]^2\  \overline{\uD^\alpha_\eta \hat f_k(\eta)}\ \hat u_l^j \ \uD^{\alpha_1}_\eta(\eta_j)\ \uD^{\alpha_2}_\eta \big(\hat f_{k-l}(\eta)\big)\  \ud \eta  \\
&=: \hat G_1 + \hat G_2,
\end{align*}
with 
$$ \hat G_1 := \Imm \sum_{\alpha,k,l} \int_{\RR^d} [A_k^\sigma(\eta)]^2\  \overline{\uD^\alpha_\eta \hat f_k(\eta)}\ \hat u_l \cdot \eta \ \uD^{\alpha}_\eta \big(\hat f_{k-l}(\eta)\big)\  \ud \eta $$
and 
$$ \hat G_2 := \Imm \sum_{j,\alpha,k,l} \int_{\RR^d} [A_k^\sigma(\eta)]^2\  \overline{\uD^\alpha_\eta \hat f_k(\eta)}\ \hat u_l^j \ \uD^{\alpha-e_j}_\eta \big(\hat f_{k-l}(\eta)\big)\  \ud \eta .$$
Now, since $u$ does not depend on $v$, we have in the physical variables
$$ \iint_{\mathbb{T}^d\times\RR^d} u(t,x) \cdot \na_v \big(A^\sigma v^\alpha f\big)^2 \ud x \ud v = 0 , $$
by Plancherel:
$$ \sum_{k,l \in \ZZ^d} \int_{\RR^d} A_k^\sigma(\eta)\ \overline{\uD^\alpha_\eta \hat f_k(\eta)} \ \hat u_l \cdot \eta \ A_{k-l}^\sigma(\eta) \ \uD^\alpha_\eta \hat f_{k-l}(\eta) \ud \eta = 0 .$$
Thus, $\hat G_1$ can be written as follow:
\begin{equation}\label{hat G_1}
\hat G_1 = \Imm \sum_{k,l \in \ZZ^d} \int_{\RR^d} A_k^\sigma(\eta)\overline{\uD^\alpha_\eta \hat f_k(\eta)} \big[A_k^\sigma(\eta)-A_{k-l}^\sigma(\eta) \big] \hat u_l \cdot \eta  \uD^\alpha_\eta \hat f_{k-l}(\eta) \ud \eta  .
\end{equation}
To estimate $\big[A_k^\sigma(\eta)-A_{k-l}^\sigma(\eta) \big]$, we will use inequalities of Lemma \ref{inegalites triangulaires}. First, we have:
$$ A_k^\sigma(\eta)-A_{k-l}^\sigma(\eta) = \big[\langle k,\eta\rangle^{\sigma}-\langle k-l,\eta\rangle^{\sigma} \big] \ue^{\lambda \langle k-l,\eta\rangle^{s}} + \big[\ue^{\lambda \langle k,\eta\rangle^{s}}-\ue^{\lambda \langle k-l,\eta\rangle^{s}}\big] \langle k,\eta\rangle^{\sigma} =: I + J . $$
Let's start by estimating $\langle k,\eta\rangle^{\sigma}-\langle k-l,\eta\rangle^{\sigma}$. By applying the mean value theorem to the function $t\mapsto t^{\frac{\sigma}{2}}$ between $X=\langle k,\eta\rangle^{2}$ and $Y=\langle k-l,\eta\rangle^{2}$ then, there exists $\theta:=\theta_{k,l,\eta} \in ]0,1[$ such that:
\begin{align*}
\langle k,\eta\rangle^{\sigma}-\langle k-l,\eta\rangle^{\sigma} &= \frac{\sigma}{2} \big[ \theta \langle k,\eta\rangle^{2}+ (1-\theta)\langle k-l,\eta\rangle^{2}\big]^{\frac{\sigma}{2}-1}\big(\langle k,\eta\rangle^{2}-\langle k-l,\eta\rangle^{2}\big)  \\
&= \frac{\sigma}{2} \bigg[ \big(\theta \langle k,\eta\rangle^{2}+ (1-\theta)\langle k-l,\eta\rangle^{2}\big)^{\frac{\sigma}{2}-1}- \langle k-l,\eta \rangle^{\sigma-2}\bigg]\big(\langle k,\eta\rangle^{2}-\langle k-l,\eta\rangle^{2}\big)  \\
&\quad + \frac{\sigma}{2}\langle k-l,\eta\rangle^{\sigma-2} \big(\langle k,\eta\rangle^{2}-\langle k-l,\eta\rangle^{2}\big)  \\
&= B + \sigma \langle k-l,\eta\rangle^{\sigma-2} \sum_{j=1}^d l_j (k_j-l_j), 
\end{align*}
where
$$ B := \frac{\sigma}{2} \bigg[ \big(\theta \langle k,\eta\rangle^{2}+ (1-\theta)\langle k-l,\eta\rangle^{2}\big)^{\frac{\sigma}{2}-1}- \langle k-l,\eta \rangle^{\sigma-2}\bigg]\big(\langle k,\eta\rangle^{2}-\langle k-l,\eta\rangle^{2}\big) + \frac{\sigma}{2}\langle k-l,\eta\rangle^{\sigma-2}|l|^2 . $$
Applying the mean value theorem once again, this time to the function $t\mapsto t^{\frac{\sigma}{2}-1}$ between $X'=\theta \langle k,\eta\rangle^{2}+ (1-\theta)\langle k-l,\eta\rangle^{2}$ and $Y'=\langle k-l,\eta\rangle^{2}$ then, there exists $\theta':=\theta_{k,l,\eta,\theta}' \in ]0,1[$ such that:
\begin{align*}
&\big(\theta \langle k,\eta\rangle^{2}+ (1-\theta)\langle k-l,\eta\rangle^{2}\big)^{\frac{\sigma}{2}-1}- \langle k-l,\eta \rangle^{\sigma-2} \\
&= \big(\frac{\sigma}{2}-1\big)\theta \big(\langle k,\eta\rangle^{2}-\langle k-l,\eta\rangle^{2}\big)\bigg[\theta \theta'\langle k,\eta\rangle^{2}+(1-\theta\theta')\langle k-l,\eta\rangle^{2}\bigg]^{\frac{\sigma}{2}-2} .
\end{align*}
By inequality 3 of Lemma \ref{inegalites triangulaires} we write:
$$\big|\langle k,\eta\rangle^{2}-\langle k-l,\eta\rangle^{2}\big| \leqslant \langle l \rangle \big(\langle k,\eta\rangle +\langle k-l,\eta\rangle\big). $$
Then, by inequality 6 of Lemma \ref{inegalites triangulaires} we obtain:
$$ \bigg|\theta \theta'\langle k,\eta\rangle^{2}+(1-\theta\theta')\langle k-l,\eta\rangle^{2}\bigg|^{\frac{\sigma}{2}-2} \lesssim \langle k,\eta\rangle^{\sigma-4} + \langle k-l,\eta\rangle^{\sigma-4} . $$
Therefore,
\begin{align*}
|B| &\lesssim \langle l \rangle^2 \big(\langle k,\eta\rangle^{2} +\langle k-l,\eta\rangle^{2}\big)\big(\langle k,\eta\rangle^{\sigma-4} +\langle k-l,\eta\rangle^{\sigma-4}\big) +\langle l \rangle^2 \langle k-l,\eta\rangle^{\sigma-2}  \\
&\lesssim \langle l \rangle^2 \big(\langle k,\eta\rangle^{\sigma-2} + \langle k,\eta\rangle^{\sigma-4} \langle k-l,\eta\rangle^{2} + \langle k,\eta\rangle^{\sigma-4}\langle k-l,\eta\rangle^{2}+\langle k,\eta\rangle^{\sigma-2}\big) .
\end{align*}
Finally, by using inequality 2 of Lemma \ref{inegalites triangulaires}, $\langle k,\eta\rangle^{2} \lesssim \langle l \rangle^2 + \langle k-l,\eta\rangle^{2}$, we get:
\begin{equation}\label{estmation de B}
|B| \lesssim \langle l \rangle^2 \big(\langle l\rangle^{\sigma-2} +\langle k-l,\eta\rangle^{\sigma-2}\big).
\end{equation}
In summary, 
\begin{equation}\label{I}
I = \bigg[ B + \sigma \langle k-l,\eta\rangle^{\sigma-2} \sum_{j=1}^d l_j (k_j-l_j) \bigg] \ue^{\lambda \langle k-l,\eta\rangle^{s}} ,
\end{equation}
with $B$ satisfying inequality \eqref{estmation de B}. \\
Now, for $J$, we write:
$$ J := \big[ \ue^{\lambda \langle k,\eta\rangle^{s}}-\ue^{\lambda \langle k-l,\eta\rangle^{s}} \big] \langle k,\eta\rangle^{\sigma} = \ue^{\lambda \langle k-l,\eta\rangle^{s}}\big[ \ue^{\lambda [\langle k,\eta\rangle^{s}-\langle k-l,\eta\rangle^{s}]}-1 \big] \langle k,\eta\rangle^{\sigma} .$$
Since $|\ue^x-1|\leqslant |x| \ue^{|x|}$ for all $x\in \RR$ then, by using inequality 5 of Lemma \ref{inegalites triangulaires}, we get
$$ |J| \leqslant \lambda \big|\langle k,\eta\rangle^{s}-\langle k-l,\eta\rangle^{s}\big| \ue^{\lambda \langle l\rangle^{s}}\ue^{\lambda \langle k-l,\eta\rangle^{s}}  \langle k,\eta\rangle^{\sigma} .$$
Now, by using inequality 7 of Lemma \ref{inegalites triangulaires} we obtain
$$ \big|\langle k,\eta\rangle^{s}-\langle k-l,\eta\rangle^{s}\big| \lesssim   \langle l \rangle \langle k,\eta\rangle^{s-1} , $$
and by inequality 6 of the same Lemma \ref{inegalites triangulaires},
$$  \langle k,\eta\rangle^{\sigma+\frac{s}{2}-1} \lesssim \langle l \rangle^{\sigma+\frac{s}{2}-1} + \langle k-l,\eta \rangle^{\sigma+\frac{s}{2}-1} . $$
Then, thanks to the last two inequalities, we write
\begin{align*}
\big|\langle k,\eta\rangle^{s}-\langle k-l,\eta\rangle^{s}\big|\langle k,\eta\rangle^{\sigma} &\lesssim \langle l \rangle \langle k,\eta\rangle^{s-1} \langle k,\eta\rangle^{\sigma}  \\
&= \langle k,\eta\rangle^{\frac{s}{2}} \langle k,\eta\rangle^{\frac{s}{2}+\sigma-1} \langle l \rangle \\
&\lesssim  \langle k,\eta\rangle^{\frac{s}{2}} \big[\langle l \rangle^{\sigma+\frac{s}{2}} +\langle l \rangle \langle k-l,\eta \rangle^{\frac{s}{2}+\sigma-1} \big].
\end{align*}
Hence,
\begin{equation}\label{J}
|J| \lesssim \lambda\langle k,\eta\rangle^{\frac{s}{2}} \big[\langle l \rangle^{\sigma+\frac{s}{2}} +\langle l \rangle \langle k-l,\eta \rangle^{\frac{s}{2}+\sigma-1} \big] \ue^{\lambda \langle l\rangle^{s}}\ue^{\lambda \langle k-l,\eta\rangle^{s}} .
\end{equation} 

\ni \textbf{Estimation of $\hat G_1$.}  Since $ A_k^\sigma(\eta)-A_{k-l}^\sigma(\eta) =I+J$ and thanks to equality \eqref{I}, we can decompose $\hat G_1$ as follows:
$$\hat G_1 = \hat G_{I,1} + \hat G_{I,2} + \hat G_J ,$$
where 
$$ \hat G_{I,1} :=  \Imm  \sum_{\alpha,k,l} \int_{\RR^d} \langle k,\eta\rangle^{\sigma} \ue^{\lambda \langle k,\eta\rangle^{s}} \overline{\uD^\alpha_\eta \hat f_k(\eta)} \  \hat u_l \cdot \eta \ B \ \ue^{\lambda \langle k-l,\eta\rangle^{s}}  \uD^\alpha_\eta \hat f_{k-l}(\eta) \ \ud \eta , $$
$$ \hat G_{I,2} := \sigma \Imm \sum_{\alpha,k,l} \sum_{j=1}^d \int_{\RR^d} \langle k,\eta\rangle^{\sigma}\ue^{\lambda \langle k,\eta\rangle^{s}} \overline{\uD^\alpha_\eta \hat f_k(\eta)} \ l_j \ \hat u_l \cdot \eta (k_j-l_j) \langle k-l,\eta\rangle^{\sigma-2}\ue^{\lambda \langle k-l,\eta\rangle^{s}} \uD^\alpha_\eta \hat f_{k-l}(\eta) \ \ud \eta $$
and 
$$ \hat G_{J} := \Imm \sum_{\alpha,k,l} \int_{\RR^d} J \ A_k^\sigma(\eta)\ \overline{\uD^\alpha_\eta \hat f_k(\eta)} \ \hat u_l \cdot \eta \ \uD^\alpha_\eta \hat f_{k-l}(\eta) \ \ud \eta .$$

\ni \textbf{Estimation of $\hat G_{I,1}$.} We have by inequality \eqref{estmation de B}
$$ \big| \hat G_{I,1} \big| \lesssim   \sum_{\alpha,k,l} \int_{\RR^d} \langle k,\eta\rangle^{\sigma} \ue^{\lambda \langle k,\eta\rangle^{s}}\big| \uD^\alpha_\eta \hat f_k(\eta)\big| \langle l\rangle^2 |\hat u_l| \big[\langle l\rangle^{\sigma-2}+\langle k-l,\eta\rangle^{\sigma-2}\big] |\eta| \ue^{\lambda \langle k-l,\eta\rangle^{s}} \big|\uD^\alpha_\eta \hat f_{k-l}(\eta) \big| \ud \eta . $$
Then, by using the inequality $\ue^x \leqslant \ue + x^2 \ue^x $ for all $x\geqslant 0$, we get
$$\big| \hat G_{I,1} \big| \lesssim  \big| \hat G_{I,11} \big| +\big| \hat G_{I,12} \big| +\big| \hat G_{I,13} \big| + \big| \hat G_{I,14} \big|, $$
where
\begin{align*}
\big| \hat G_{I,11} \big| &\lesssim   \sum_{\alpha,k,l} \int_{\RR^d} \langle k,\eta\rangle^{\sigma} \ue^{\lambda \langle k,\eta\rangle^{s}}\big| \uD^\alpha_\eta \hat f_k(\eta)\big| \langle l\rangle^{\sigma} |\hat u_l| |\eta| \big|\uD^\alpha_\eta \hat f_{k-l}(\eta) \big| \ud \eta , \\
\big| \hat G_{I,12} \big| &\lesssim  \lambda^2 \sum_{\alpha,k,l} \int_{\RR^d} \langle k,\eta\rangle^{\sigma} \ue^{\lambda \langle k,\eta\rangle^{s}}\big| \uD^\alpha_\eta \hat f_k(\eta)\big| \langle l\rangle^{\sigma} |\hat u_l| |\eta| \langle k-l,\eta\rangle^{2s} \ue^{\lambda \langle k-l,\eta\rangle^{s}}\big|\uD^\alpha_\eta \hat f_{k-l}(\eta) \big| \ud \eta ,\\
\big| \hat G_{I,13} \big| &\lesssim   \sum_{\alpha,k,l} \int_{\RR^d} \langle k,\eta\rangle^{\sigma} \ue^{\lambda \langle k,\eta\rangle^{s}}\big| \uD^\alpha_\eta \hat f_k(\eta)\big| \langle l\rangle^{2} |\hat u_l| |\eta| \langle k-l,\eta\rangle^{\sigma-2} \big|\uD^\alpha_\eta \hat f_{k-l}(\eta) \big| \ud \eta ,\\
 \big| \hat G_{I,14} \big| &\lesssim \lambda^2  \sum_{\alpha,k,l} \int_{\RR^d} \langle k,\eta\rangle^{\sigma} \ue^{\lambda \langle k,\eta\rangle^{s}}\big| \uD^\alpha_\eta \hat f_k(\eta)\big| \langle l\rangle^{2} |\hat u_l| |\eta| \langle k-l,\eta\rangle^{\sigma+2s-2} \ue^{\lambda \langle k-l,\eta\rangle^{s}} \big|\uD^\alpha_\eta \hat f_{k-l}(\eta) \big| \ud \eta .
\end{align*}
Note that $|\eta|\leqslant \langle k-l,\eta\rangle$. Now, by applying inequality \eqref{Young 2} of Lemma \ref{inegalite de Young} for the first two inequalities, with $\gamma=1$ and $\beta=0$ for $\hat G_{I,11}$, and  with $\gamma=2s+1$ and $\beta=\frac{s}{2}$ for $\hat G_{I,12}$, for $\sigma > \frac{d}{2}+1$,  we obtain
\begin{equation}\label{G_I11}
\big| \hat G_{I,11} \big| \lesssim \| u\|_\sigma \| f \|_{\sigma,M} \|f\|_{\lambda,\sigma,M,s} ,
\end{equation}
and for $\sigma > \frac{d}{2}+\frac{3s}{2}+1$,  we obtain
\begin{equation}\label{G_I12}
\big| \hat G_{I,12} \big| \lesssim \lambda^2 \| u\|_\sigma \|f\|_{\lambda,\sigma,M,s} \|f\|_{\lambda,\sigma+\frac{s}{2},M,s} .
\end{equation}
For $\hat G_{I,13}$, by applying inequality \eqref{Young 1} for $\nu=2$ and $\beta=-1$, for $\sigma > \frac{d}{2}+2$ we get
\begin{equation}\label{G_I13}
\big| \hat G_{I,13} \big| \lesssim \| u\|_\sigma \| f \|_{\sigma-1,M} \|f\|_{\lambda,\sigma,M,s} .
\end{equation}
For $\hat G_{I,14}$, since $\langle k-l,\eta\rangle^{\frac{s}{2}} \leqslant \langle k,\eta\rangle^{\frac{s}{2}} + \langle l\rangle^{\frac{s}{2}}$ (by inequality 5 of Lemma \ref{inegalites triangulaires}) and $s\in (0,1]$ then,
\begin{align*}
\big| \hat G_{I,14} \big| &\lesssim \lambda^2  \sum_{\alpha,k,l} \int_{\RR^d} \langle k,\eta\rangle^{\sigma+\frac{s}{2}} \ue^{\lambda \langle k,\eta\rangle^{s}}\big| \uD^\alpha_\eta \hat f_k(\eta)\big| \langle l\rangle^{2} |\hat u_l| \langle k-l,\eta\rangle^{\sigma+\frac{s}{2}} \ue^{\lambda \langle k-l,\eta\rangle^{s}} \big|\uD^\alpha_\eta \hat f_{k-l}(\eta) \big| \ud \eta \\
&+ \lambda^2  \sum_{\alpha,k,l} \int_{\RR^d} \langle k,\eta\rangle^{\sigma} \ue^{\lambda \langle k,\eta\rangle^{s}}\big| \uD^\alpha_\eta \hat f_k(\eta)\big| \langle l\rangle^{2+\frac{s}{2}} |\hat u_l| |\eta| \langle k-l,\eta\rangle^{\sigma+\frac{s}{2}} \ue^{\lambda \langle k-l,\eta\rangle^{s}} \big|\uD^\alpha_\eta \hat f_{k-l}(\eta) \big| \ud \eta .
\end{align*}
Thus, by applying inequality \eqref{Young 2} to the first line in the previous inequality, for $\nu=2$ and $\beta=\frac{s}{2}$, and by applying inequality \eqref{Young 1} to the second line, for $\nu=2+\frac{s}{2}$ and $\beta=\frac{s}{2}$, we get for $\sigma > \frac{d}{2}+\frac{s}{2}+2$ 
$$ \big| \hat G_{I,14} \big| \lesssim  \lambda^2 \left( \| u\|_\sigma \| f \|_{\lambda,\sigma,M,s} \|f\|_{\lambda,\sigma+\frac{s}{2},M,s} + \| u\|_{2+\frac{s}{2}} \|f\|_{\lambda,\sigma+\frac{s}{2},M,s}^2 \right) .$$
Hence, for $\sigma > \frac{d}{2}+\frac{s}{2}+2$
\begin{equation}\label{G_I14}
\big| \hat G_{I,14} \big| \lesssim \lambda^2 \| u\|_\sigma \|f\|_{\lambda,\sigma+\frac{s}{2},M,s}^2.
\end{equation}
Therefore, by summing inequalities \eqref{G_I11}, \eqref{G_I12}, \eqref{G_I13} and \eqref{G_I14} we obtain for $\sigma > \frac{d}{2}+\frac{s}{2}+2$:
\begin{equation}\label{G_I1}
\big| \hat G_{I,1} \big| \lesssim  \| u\|_\sigma \left(\|f\|_{\sigma,M}\|f\|_{\lambda,\sigma,M,s}+ \lambda^2 \|f\|_{\lambda,\sigma+\frac{s}{2},M,s}^2 \right).
\end{equation}

\ni \textbf{Estimation of $\hat G_{I,2}$.} Recall that $\hat G_{I,2}$ is defined by
$$ \hat G_{I,2} := \sigma \Imm \sum_{\alpha,k,l} \sum_{j=1}^d \int_{\RR^d} \langle k,\eta\rangle^{\sigma}\ue^{\lambda \langle k,\eta\rangle^{s}} \overline{\uD^\alpha_\eta \hat f_k(\eta)} \ l_j \ \hat u_l \cdot \eta (k_j-l_j) \langle k-l,\eta\rangle^{\sigma-2}\ue^{\lambda \langle k-l,\eta\rangle^{s}} \uD^\alpha_\eta \hat f_{k-l}(\eta) \ \ud \eta .$$
Thus, $\hat G_{I,2}$ can be seen, using the inverse Fourier transform, as
$$ \hat G_{I,2} = \sigma \Imm \sum_{\alpha,k} \sum_{j=1}^d \int_{\RR^d} \overline{\mathcal{F}\big( A^\sigma v^\alpha f\big)} \mathcal{F}\bigg( \frac{\pa u}{\pa x_j} \cdot \na_v \frac{\pa}{\pa x_j}\big(A^{\sigma-2} v^\alpha f\big)\bigg) \ \ud \eta . $$
Hence,
\begin{equation}\label{G_I2}
\big| \hat G_{I,2} \big| \lesssim  \| \na_x u\|_\infty \|f\|_{\lambda,\sigma,M,s}^2 .
\end{equation}

Note that each time we want to have $L^\infty$ estimates (for $u$), we go back through Fourier using Parseval.  \\

\ni \textbf{Estimation of $\hat G_{J}$.} By using inequality \eqref{J},  we write  
\begin{align*}
\big| \hat G_{J} \big| &\leqslant \sum_{\alpha,k,l} \int_{\RR^d} |J| A_k^\sigma(\eta)\big|\uD^\alpha_\eta \hat f_k(\eta)\big| |\hat u_l| |\eta| \big|\uD^\alpha_\eta \hat f_{k-l}(\eta)\big| \ud \eta \\
&\lesssim \lambda \sum_{\alpha,k,l} \int_{\RR^d} \langle k,\eta\rangle^{\sigma+\frac{s}{2}}\ue^{\lambda \langle k,\eta\rangle^{s}}\big|\uD^\alpha_\eta \hat f_k(\eta)\big| \langle l\rangle^{\sigma+\frac{s}{2}}\ue^{\lambda \langle l\rangle^{s}}|\hat u_l| |\eta| \ue^{\lambda \langle k-l,\eta\rangle^{s}} \big|\uD^\alpha_\eta \hat f_{k-l}(\eta)\big| \ud \eta \\
&+\lambda \sum_{\alpha,k,l} \int_{\RR^d} \langle k,\eta\rangle^{\sigma+\frac{s}{2}}\ue^{\lambda \langle k,\eta\rangle^{s}}\big|\uD^\alpha_\eta \hat f_k(\eta)\big| \langle l\rangle \ue^{\lambda \langle l\rangle^{s}}|\hat u_l| \langle k-l,\eta\rangle^{\sigma+\frac{s}{2}}\ue^{\lambda \langle k-l,\eta\rangle^{s}} \big|\uD^\alpha_\eta \hat f_{k-l}(\eta)\big| \ud \eta .
\end{align*}
We proceed as in $\hat G_{I,1}$, using the inequality $\ue^x \leqslant 1+x\ue^x$ for all $x\geqslant0$, we get:
$$ \big| \hat G_{J} \big| \lesssim \big| \hat G_{J1} \big| + \big| \hat G_{J2} \big| + \big| \hat G_{J3} \big| + \big| \hat G_{J4} \big| ,$$
where
\begin{align*}
\big| \hat G_{J1} \big| &\lesssim  \lambda \sum_{\alpha,k,l} \int_{\RR^d} \langle k,\eta\rangle^{\sigma+\frac{s}{2}} \ue^{\lambda \langle k,\eta\rangle^{s}} \big|\uD^\alpha_\eta \hat f_k(\eta)\big| \langle l\rangle^{\sigma+\frac{s}{2}} \ue^{\lambda \langle l\rangle^{s}}|\hat u_l| |\eta| \big|\uD^\alpha_\eta \hat f_{k-l}(\eta)\big| \ud \eta ,\\
\big| \hat G_{J2} \big| &\lesssim  \lambda^2 \sum_{\alpha,k,l} \int_{\RR^d} \langle k,\eta\rangle^{\sigma+\frac{s}{2}} \ue^{\lambda \langle k,\eta\rangle^{s}}\big|\uD^\alpha_\eta \hat f_k(\eta)\big| \langle l\rangle^{\sigma+\frac{s}{2}} \ue^{\lambda \langle l\rangle^{s}}|\hat u_l| \langle k-l,\eta\rangle^{1+s} \ue^{\lambda \langle k-l,\eta\rangle^{s}} \big|\uD^\alpha_\eta \hat f_{k-l}(\eta)\big| \ud \eta ,\\ 
\big| \hat G_{J3} \big| &\lesssim  \lambda \sum_{\alpha,k,l} \int_{\RR^d} \langle k,\eta\rangle^{\sigma+\frac{s}{2}} \ue^{\lambda \langle k,\eta\rangle^{s}}\big|\uD^\alpha_\eta \hat f_k(\eta)\big| \langle l\rangle |\hat u_l|  \langle k-l,\eta\rangle^{\sigma+\frac{s}{2}} \big|\uD^\alpha_\eta \hat f_{k-l}(\eta)\big| \ud \eta ,  \\ 
\big| \hat G_{J4} \big| &\lesssim \lambda^2 \sum_{\alpha,k,l} \int_{\RR^d} \langle k,\eta\rangle^{\sigma+\frac{s}{2}} \ue^{\lambda \langle k,\eta\rangle^{s}}\big|\uD^\alpha_\eta \hat f_k(\eta)\big| \langle l\rangle^{1+s} \ue^{\lambda \langle l\rangle^{s}}|\hat u_l|  \langle k-l,\eta\rangle^{\sigma+\frac{s}{2}} \ue^{\lambda \langle k-l,\eta\rangle^{s}} \big|\uD^\alpha_\eta \hat f_{k-l}(\eta)\big| \ud \eta  .
\end{align*} 
By applying inequality \eqref{Young 2} to the first two inequalities, we obtain for $\sigma > \frac{d}{2}+1$,
\begin{equation}\label{G_J1}
\big| \hat G_{J1} \big| \lesssim \lambda \|f\|_{\sigma,M} \|u\|_{\lambda,\sigma+\frac{s}{2},s} \|f\|_{\lambda,\sigma+\frac{s}{2},M,s} ,
\end{equation}
and for $\sigma > \frac{d}{2}+s+1$,
\begin{equation}\label{G_J2}
\big| \hat G_{J2} \big| \lesssim \lambda^2 \|f\|_{\lambda,\sigma,M,s} \|u\|_{\lambda,\sigma+\frac{s}{2},s} \|f\|_{\lambda,\sigma+\frac{s}{2},M,s}  .
\end{equation}
For $\hat G_{J3}$ and $\hat G_{I,14}$, we have thanks to inequality \eqref{Young 1}, for $\sigma > \frac{d}{2}+1$
\begin{equation}\label{G_J3}
\big| \hat G_{J3} \big| \lesssim \lambda \|u\|_{\sigma} \|f\|_{\lambda,\sigma+\frac{s}{2},M,s}^2 ,
\end{equation}
and for $\sigma > \frac{d}{2}+s+1$, 
\begin{equation}\label{G_J4}
\big| \hat G_{J4} \big| \lesssim \lambda^2 \|u\|_{\lambda,\sigma,s} \|f\|_{\lambda,\sigma+\frac{s}{2},M,s}^2  .
\end{equation}
Thus, by summing the last four inequalities, we get:
\begin{align}\label{G_J}
\big| \hat G_{J} \big| &\lesssim \left(\lambda \|u\|_{\sigma}+\lambda^2 \|u\|_{\lambda,\sigma,s}\right) \|f\|_{\lambda,\sigma+\frac{s}{2},M,s}^2  \nonumber \\
&+  \left(\lambda \|f\|_{\sigma,M} + \lambda^2 \|f\|_{\lambda,\sigma,M,s}\right) \|u\|_{\lambda,\sigma+\frac{s}{2},s} \|f\|_{\lambda,\sigma+\frac{s}{2},M,s} .
\end{align} 
Finally, by summing \eqref{G_I1}, \eqref{G_I2} and \eqref{G_J}, we obtain the following estimate for $\hat G_1$:
\begin{align}\label{G_1}
\big| \hat G_1 \big| &\lesssim   \| \na_x u\|_\infty \|f\|_{\lambda,\sigma,M,s}^2 + \|u\|_{\sigma} \|f\|_{\sigma,M} \|f\|_{\lambda,\sigma,M,s} \nonumber \\
&+\left( (\lambda+\lambda^2) \|u\|_{\sigma} + \lambda^2 \|u\|_{\lambda,\sigma,s}\right) \|f\|_{\lambda,\sigma+\frac{s}{2},M,s}^2  \nonumber \\
&+  \left( \lambda \|f\|_{\sigma,M} + \lambda^2 \|f\|_{\lambda,\sigma,M,s}\right) \|u\|_{\lambda,\sigma+\frac{s}{2},s} \|f\|_{\lambda,\sigma+\frac{s}{2},M,s} .
\end{align}
Thus, it only remains to estimate $\hat G_2$ to conclude the Gevrey estimates for the solution of the Vlasov equation. \\

\ni \textbf{Estimation of $\hat G_2$.} Recall that $\hat G_2$ is defined by 
$$ \hat G_2 := \Imm \sum_{\alpha,k,l} \sum_{j=1}^d \int_{\RR^d} [A_k^\sigma(\eta)]^2\  \overline{\uD^\alpha_\eta \hat f_k(\eta)}\ \hat u_l^j \ \uD^{\alpha-e_j}_\eta \big(\hat f_{k-l}(\eta)\big)\  \ud \eta ,$$
which we decompose as follows:
\begin{align*}
\hat G_2 &= \Imm \sum_{j,\alpha,k,l}  \int_{\RR^d} A_k^\sigma(\eta)\  \overline{\uD^\alpha_\eta \hat f_k(\eta)}\ \hat u_l^j \ A_{k-l}^\sigma(\eta) \uD^{\alpha-e_j}_\eta \big(\hat f_{k-l}(\eta)\big)\  \ud \eta \\
& + \Imm \sum_{j,\alpha,k,l} \int_{\RR^d} A_k^\sigma(\eta)\  \overline{\uD^\alpha_\eta \hat f_k(\eta)}\ \hat u_l^j \big[A_k^\sigma(\eta)-A_{k-l}^\sigma(\eta)\big] \uD^{\alpha-e_j}_\eta \big(\hat f_{k-l}(\eta)\big)\  \ud \eta \\
&=: \hat G_{21} + \hat G_{22} .
\end{align*}
For $\hat G_{21}$, by the inverse Fourier transform theorem, we write:
\begin{align}\label{G_21}
\big|\hat G_{21} \big| &:= \bigg| \Imm  \sum_{j,\alpha,k,l}  \int_{\RR^d} A_k^\sigma(\eta)\  \overline{\uD^\alpha_\eta \hat f_k(\eta)}\ \hat u_l^j \ A_{k-l}^\sigma(\eta) \uD^{\alpha-e_j}_\eta \big(\hat f_{k-l}(\eta)\big)\  \ud \eta \bigg| \nonumber \\
&= \bigg| \Imm \sum_{j,\alpha,k} \int_{\RR^d} \overline{\mathcal{F} \left(A^\sigma v^\alpha f\right)}\ \mathcal{F} \left(u^j A^{\sigma} v^{\alpha-e_j}  f\right)\  \ud \eta \bigg| \nonumber \\
&\lesssim  \|u\|_\infty \|f\|_{\lambda,\sigma,M,s}^2 .
\end{align}
For the term $\hat G_{22}$, proceeding exactly in the same way as for $\hat G_1$, we obtain for $\sigma > \frac{d}{2}+2$:
\begin{align}\label{G_22}
\big| \hat G_{22} \big| &\lesssim   \| \na_x u\|_\infty \|f\|_{\lambda,\sigma,M,s}^2 + \|u\|_{\sigma} \|f\|_{\sigma,M} \|f\|_{\lambda,\sigma,M,s} \nonumber \\
&+\left( (\lambda+\lambda^2) \|u\|_{\sigma} + \lambda^2 \|u\|_{\lambda,\sigma,s}\right) \|f\|_{\lambda,\sigma,M,s}^2  \nonumber \\
&+  \left( \lambda \|f\|_{\sigma,M} + \lambda^2 \|f\|_{\lambda,\sigma,M,s}\right) \|u\|_{\lambda,\sigma+\frac{s}{2},s} \|f\|_{\lambda,\sigma+\frac{s}{2},M,s} .
\end{align}
Finally, returning to \eqref{d/dt |f|^2 = lambda'|f|^2+E+F+G} and summing all the inequalities obtained on $\hat E$, $\hat F$ and $\hat G$, namely: \eqref{hat E}, \eqref{hat F}, \eqref{G_1}, \eqref{G_21} and \eqref{G_22}, we obtain the following estimate for $\sigma > \frac{d}{2}+\frac{s}{2}+2$
\begin{align}
\frac{1}{2}\frac{\ud}{\ud t} \|f\|_{\lambda,\sigma,M,s}^2 &\lesssim \left( \|u\|_{W^{1,\infty}}+1\right) \|f\|_{\lambda,\sigma,M,s}^2 + \|u\|_{\sigma} \|f\|_{\sigma,M} \|f\|_{\lambda,\sigma,M,s}  \nonumber \\
&+  \left( \dot{\lambda} + \lambda (1+ \|u\|_{\sigma}) + \lambda^2 \|u\|_{\lambda,\sigma,s}\right) \|f\|_{\lambda,\sigma+\frac{s}{2},M,s}^2  \nonumber \\
&+ \left( \lambda \|f\|_{\sigma,M} + \lambda^2 \|f\|_{\lambda,\sigma,M,s} \right) \|u\|_{\lambda,\sigma+\frac{s}{2},s} \|f\|_{\lambda,\sigma+\frac{s}{2},M,s} .
\end{align}
This completes the proof of Proposition \ref{Gevrey Vlasov}.
\ep
\subsection{Gevrey estimates for the Navier-Stokes field}
The purpose of this subsection is to prove the inequality of the following proposition:
\begin{proposition}\label{Gevrey NS}
Let $s\in (0,1]$, $M>\frac{d}{2}+1$ and $\sigma>\frac{d}{2}+\frac{s}{2}+2$. For $f \in \mathcal{G}^{\lambda,\sigma+\frac{s}{2},M,\frac{1}{s}}(\mathbb{T}^d\times\RR^d)$ and $u \in \mathcal{G}^{\lambda,\sigma+\frac{s}{2},\frac{1}{s}}(\mathbb{T}^d)$ satisfying \eqref{VNS}, the following estimate holds
\begin{align}\label{estimation Gevrey pour u}
\frac{1}{2}\frac{\ud}{\ud t} \| u \|_{\lambda,\sigma,s}^2 &+ \| u \|_{\lambda,\sigma+1,s}^2 \lesssim \left(\| \na_x u\|_{\infty}+\| f\|_{\sigma,M}+\lambda^2 \| f\|_{\lambda,\sigma,M,s}\right) \| u \|_{\lambda,\sigma,s}^2 + \| u \|_\sigma^2\| u \|_{\lambda,\sigma,s}  \nonumber \\
 &+\| f\|_{\lambda,\sigma,M,s} \| u \|_{\lambda,\sigma,s}+ \left( \dot{\lambda}+\lambda \| u\|_\sigma + \lambda^2 (\| u \|_\sigma+\|u\|_{\lambda,\sigma,s}) \right) \| u \|_{\lambda,\sigma+\frac{s}{2},s}^2.
\end{align}
\end{proposition}

As in Proposition \ref{Gevrey Vlasov} for Vlasov,  we will choose a suitable $\lambda$ to absorb the term $\| u \|_{\lambda,\sigma+\frac{s}{2},s}$.  More precisely, we will take $\lambda$ such that the following inequality is satisfied
\begin{equation}\label{contrainte 2 lambda}
 \dot{\lambda}+\lambda \| u\|_\sigma + \lambda^2 (\| u \|_\sigma+\|u\|_{\lambda,\sigma,s}) \leqslant 0.
\end{equation}
Since we have two constraints now \eqref{contrainte 1 lambda} and \eqref{contrainte 2 lambda}, we have to take a $\lambda$ which satisfies both at the same time, but here the same $\lambda$ satisfying \eqref{contrainte 1 lambda} is good and it gives :
$$ \frac{1}{2}\frac{\ud}{\ud t} \| u \|_{\lambda,\sigma,s}^2 + \| u \|_{\lambda,\sigma+1,s}^2 \lesssim \left(\| \na_x u\|_{\infty}+\| f\|_{\sigma,M}+\lambda^2 \| f\|_{\lambda,\sigma,M,s}\right) \| u \|_{\lambda,\sigma,s}^2 + \| u \|_\sigma^2\| u \|_{\lambda,\sigma,s} .$$
Therefore, for $\| f\|_{\lambda,\sigma,M,s}$ finite, Gronwall's lemma gives us a control of the norm $ \| u \|_{\lambda,\sigma,s}$. \\

\ni \bpp \ref{Gevrey NS}. Recall that $\ds \Lambda := (\mathrm{Id}-\Delta_x)^{\frac{1}{2}}$ and that 
$$ \|u\|_{\lambda,\sigma,s}^2 := \int_{\mathbb{T}^d} \big|\Lambda^{\sigma} \ue^{\lambda \Lambda^{s}} u \big|^2 \ud x = \sum_{k \in \ZZ^d} \langle k\rangle^{2\sigma} \ue^{2\lambda \langle k\rangle^{s}} |\hat u_k|^2 .$$
Then, applying $\Lambda^{\sigma} \ue^{\lambda \Lambda^{s}}$ to the equation satisfied by $u$ and taking the scalar product with $\Lambda^{\sigma} \ue^{\lambda \Lambda^{s}} u$, we obtain:
\begin{align*}
&\big\langle \Lambda^{\sigma} \ue^{\lambda \Lambda^{s}} \pa_t u,\Lambda^{\sigma} \ue^{\lambda \Lambda^{s}} u \big\rangle_{L^2} + \| \Lambda^{\sigma+1} \ue^{\lambda \Lambda^{s}}u \|_{L^2}^2+ \big\langle \Lambda^{\sigma} \ue^{\lambda \Lambda^{s}}(u\cdot\na_x u),\Lambda^{\sigma} \ue^{\lambda \Lambda^{s}} u \big\rangle_{L^2} \\
&= \big\langle \Lambda^{\sigma} \ue^{\lambda \Lambda^{s}} u , \Lambda^{\sigma} \ue^{\lambda \Lambda^{s}} j_f \big\rangle_{L^2} - \big\langle \Lambda^{\sigma} \ue^{\lambda \Lambda^{s}} (\rho_f u),\Lambda^{\sigma} \ue^{\lambda \Lambda^{s}} u\big\rangle_{L^2} .
\end{align*}
We have $\| \Lambda^{\sigma+1} \ue^{\lambda \Lambda^{s}}u \|_{L^2}^2 = \|u\|_{\lambda,\sigma+1,s}^2$ and 
$$\frac{1}{2} \frac{\ud}{\ud t} \|u\|_{\lambda,\sigma,s}^2 = \frac{1}{2} \frac{\ud}{\ud t} \| \Lambda^{\sigma} \ue^{\lambda \Lambda^{s}}u \|_{L^2}^2 = \dot{\lambda} \|u\|_{\lambda,\sigma+\frac{s}{2},s}^2 +\big\langle \Lambda^{\sigma} \ue^{\lambda \Lambda^{s}} \pa_t u,\Lambda^{\sigma} \ue^{\lambda \Lambda^{s}} u \big\rangle_{L^2} .$$
Therefore,
\begin{align}
\frac{1}{2} \frac{\ud}{\ud t} \|u\|_{\lambda,\sigma,s}^2 + \|u\|_{\lambda,\sigma+1,s}^2 &\leqslant  \dot{\lambda} \|u\|_{\lambda,\sigma+\frac{s}{2},s}^2 + \left| \big\langle \Lambda^{\sigma} \ue^{\lambda \Lambda^{s}}(u\cdot\na_x u),\Lambda^{\sigma} \ue^{\lambda \Lambda^{s}} u \big\rangle_{L^2} \right| \nonumber \\
&+ \|j\|_{\lambda,\sigma,s} \|u\|_{\lambda,\sigma,s} + \|\rho u \|_{\lambda,\sigma,s} \|u\|_{\lambda,\sigma,s} .
\end{align}
We claim that for $\sigma > \frac{d}{2}+\frac{s}{2}+2$, 
\begin{align}\label{estimation du commutateur}
\left| \big\langle \Lambda^{\sigma} \ue^{\lambda \Lambda^{s}}(u\cdot\na_x u),\Lambda^{\sigma} \ue^{\lambda \Lambda^{s}} u \big\rangle_{L^2} \right| &\lesssim \|\nabla_x u\|_{\infty}  \| u \|_{\lambda,\sigma,s}^2 + \| u \|_\sigma^2 \| u \|_{\lambda,\sigma,s}  \nonumber \\
 &+ \left( (\lambda+\lambda^2) \| u\|_\sigma + \lambda^2 \|u\|_{\lambda,\sigma,s} \right) \| u \|_{\lambda,\sigma+\frac{s}{2},s}^2.
\end{align}
Let denote by $H$ the term 
$$ H  :=  \big\langle \Lambda^{\sigma} \ue^{\lambda \Lambda^{s}}(u\cdot\na_x u),\Lambda^{\sigma} \ue^{\lambda \Lambda^{s}} u \big\rangle_{L^2} .$$
Since $\na_x \cdot u =0$ then,
$$ H  =  \big\langle u\cdot\na_x u,\Lambda^{2\sigma} \ue^{2\lambda \Lambda^{s}} u \big\rangle_{L^2} - \big\langle u\cdot\na_x \big(\Lambda^{\sigma} \ue^{\lambda \Lambda^{s}}u\big),\Lambda^{\sigma} \ue^{\lambda \Lambda^{s}} u \big\rangle_{L^2}.$$
Thus, in Fourier variable, we write
\begin{align*}
H  &=  \ui \sum_{k \in \ZZ^d} \bigg[\mathcal{F}[u\cdot\na_x u](k) \langle k\rangle^{2\sigma} \ue^{2\lambda \langle k\rangle^{s}} \overline{\hat u_k}  -  \mathcal{F}\big[u\cdot\na_x \big(\Lambda^{\sigma}\ue^{\lambda \Lambda^{s}} u\big)\big](k) \langle k\rangle^{\sigma} \ue^{\lambda \langle k\rangle^{s}} \overline{\hat u_k} \bigg] \\
&= \ui \sum_{k,l \in \ZZ^d} \bigg[\langle k\rangle^{\sigma} \ue^{\lambda \langle k\rangle^{s}} \overline{\hat u_k} \left( \langle k\rangle^{\sigma} \ue^{\lambda \langle k\rangle^{s}} - \langle k-l\rangle^{\sigma} \ue^{\lambda \langle k-l\rangle^{s}} \right) \hat u_l \ (k-l) \cdot \hat u_{k-l} \bigg] . 
\end{align*}
Observe that $\ds \langle k\rangle^{\sigma} \ue^{\lambda \langle k\rangle^{s}} = A_k^\sigma(0)$. Therefore,
$$ H = \ui \sum_{k,l} A_k^\sigma(0) \overline{\hat u_k} \left( A_k^\sigma(0) - A_{k-l}^\sigma(0) \right) \hat u_l \ (k-l) \cdot \hat u_{k-l} . $$
Thus, proceeding exactly the same way as for $\hat G_1$ for $\eta=0$ (see formula \eqref{hat G_1}), we obtain for $\sigma > \frac{d}{2}+\frac{s}{2}+2$:
$$
|H| \lesssim \|\nabla_x u\|_{\infty}  \| u \|_{\lambda,\sigma,s}^2 + \| u \|_\sigma^2 \| u \|_{\lambda,\sigma,s} + \left( (\lambda+\lambda^2) \| u\|_\sigma + \lambda^2 \|u\|_{\lambda,\sigma,s} \right) \| u \|_{\lambda,\sigma+\frac{s}{2},s}^2.
$$
In order to continue these estimates, we need the following lemma:
\begin{lemma}\label{commutateur, j et rho} The following estimates hold for $\sigma > \frac{d}{2}+2s$ and $M>\frac{d}{2}$, 
\begin{align}\label{estimation de rho u}
\|\rho u \|_{\lambda,\sigma,s} &\lesssim  \| u \|_{\sigma} \| f \|_{\lambda,\sigma,M,s} + \| f \|_{\sigma,M} \| u \|_{\lambda,\sigma,s} + \lambda^2 \| u \|_{\lambda,\sigma,s} \| f \|_{\lambda,\sigma,M,s}.
\end{align}
\end{lemma}

\ni \bpl \ref{commutateur, j et rho}. We have 
$$ \|\rho u \|_{\lambda,\sigma,s}^2 = \sum_{k \in \ZZ^d} \langle k\rangle^{2\sigma} \ue^{2\lambda \langle k\rangle^{s}} \big|\widehat{\rho u}_k \big|^2 =  \sum_{k \in \ZZ^d} \bigg| \sum_{l \in \ZZ^d} \langle k\rangle^{\sigma} \ue^{\lambda \langle k\rangle^{s}} \hat \rho_l \hat u_{k-l}\bigg|^2 .$$
By using inequalities 1 and 2 of Lemma \ref{inegalites triangulaires} and the inequality $\ue^x \leqslant e + x^2\ue^x$ for $x \geqslant 0$, we write:
\begin{align*}
\|\rho u \|_{\lambda,\sigma,s}^2 &\lesssim \sum_{k} \bigg| \sum_{l} [\langle l\rangle^{\sigma}+\langle k-l\rangle^{\sigma}] \ue^{\lambda \langle l\rangle^{s}}\ue^{\lambda \langle k-l\rangle^{s}} |\hat \rho_l| |\hat u_{k-l}|\bigg|^2 \\
&\lesssim \sum_{k} \bigg| \sum_{l} \langle l\rangle^{\sigma}\ue^{\lambda \langle l\rangle^{s}}|\hat \rho_l| \left(e+\lambda^2\langle k-l\rangle^{2s} \ue^{\lambda \langle k-l\rangle^{s}} \right) |\hat u_{k-l}| \\
& \hspace{1.5cm} +\left(e+\lambda^2\langle l\rangle^{2s} \ue^{\lambda \langle l\rangle^{s}} \right)|\hat \rho_l|  \langle k-l\rangle^{\sigma} \ue^{\lambda \langle k-l\rangle^{s}} |\hat u_{k-l}|\bigg|^2 .
\end{align*}
Then, by applying inequalities \eqref{Young 1} and \eqref{Young 2}, we obtain for $\sigma>\frac{d}{2}+2s$
$$ \|\rho u \|_{\lambda,\sigma,s}^2 \lesssim \|\rho \|_{\lambda,\sigma,s}^2 \| u \|_{\sigma}^2 + \|\rho\|_{\sigma}^2 \| u \|_{\lambda,\sigma,s}^2 + \lambda^2 \|\rho \|_{\lambda,\sigma,s}^2 \| u \|_{\lambda,\sigma,s}^2 .$$
Finally, we get inequality \eqref{estimation de rho u} by using  inequality \eqref{rho < f}, $\|\rho \|_{\lambda,\sigma,s}\lesssim \|f \|_{\lambda,\sigma,M,s}$ for $M>\frac{d}{2}$. 
\epl

\ni Returning now to the estimate of $\frac{1}{2} \frac{\ud}{\ud t} \| u \|_{\lambda,\sigma,s}^2$. By using Lemma \ref{commutateur, j et rho} and inequality \eqref{rho < f} then, for $\sigma > \frac{d}{2}+\frac{s}{2}+2$ and $M>\frac{d}{2}+1$, we obtain:
\begin{align*}
\frac{1}{2}\frac{\ud}{\ud t} \| u \|_{\lambda,\sigma,s}^2 + \| u \|_{\lambda,\sigma+1,s}^2 &\lesssim \left(\|\na_x u\|_{\infty}+\| f\|_{\sigma,M}+\lambda^2\| f\|_{\lambda,\sigma,M,s}\right) \| u \|_{\lambda,\sigma,s}^2 \nonumber \\
&+ \| u \|_\sigma^2\| u \|_{\lambda,\sigma,s}  \nonumber + \| f\|_{\lambda,\sigma,M,s} \| u \|_{\lambda,\sigma,s} \\
&+ \left( \dot{\lambda}+(\lambda+\lambda^2) \| u\|_\sigma + \lambda^2 \|u\|_{\lambda,\sigma,s} \right) \| u \|_{\lambda,\sigma+\frac{s}{2},s}^2 .
\end{align*}
\ep
\begin{remark}We used the inequalities $\ue^x \leqslant 1 + x \ue^x $ and $\ue^x \leqslant \ue + x^2 \ue^x $ for $x \geqslant 0$ in order to reduce the power of the Gevrey norms and apply Gronwall instead of having a power greater than $2$, which gives results for a time that depends on the initial data. More precisely, without using the last two inequalities we will have terms like $\|u\|_{\lambda,\sigma,s}\|f\|_{\lambda,\sigma,M,s}^2$ and an inequality of the type 
$$ \frac{1}{2}\frac{\ud}{\ud t} \big(\|u\|_{\lambda,\sigma,s}^2+\|f\|_{\lambda,\sigma,M,s}^2 \big) \lesssim \phi(t,\|u\|_{\sigma},\|f\|_{\sigma,M}) \big(\|u\|_{\lambda,\sigma,s}^2 + \|f\|_{\lambda,\sigma,M,s}^2 \big)^{3/2} .$$
\end{remark}

\subsection{Proof of the main results}\label{subsection main results}
\textbf{Proof of Theorem \ref{propagation Gevrey}. } 
 By Proposition \ref{Gevrey Vlasov} and Proposition \ref{Gevrey NS}, we have
\begin{align*}
&\frac{1}{2}\frac{\ud}{\ud t} \|f\|_{\lambda,\sigma,M,s}^2 \lesssim \|f\|_{\lambda,\sigma,M,s} \bigg( \left(\|u\|_{W^{1,\infty}}+1\right) \|f\|_{\lambda,\sigma,M,s} + \|u\|_{\sigma} \|f\|_{\sigma,M} \bigg) \\
&  +  \left( \dot{\lambda} + \lambda (1+ \|u\|_{\sigma}+\|f\|_{\sigma,M}) + \lambda^2 (\|u\|_{\lambda,\sigma,s}+\|f\|_{\lambda,\sigma,M,s})\right) \left(\|f\|_{\lambda,\sigma+\frac{s}{2},M,s}^2  + \|u\|_{\lambda,\sigma+\frac{s}{2},s}^2 \right) 
\end{align*}
and 
\begin{align*}
\frac{1}{2}\frac{\ud}{\ud t} \| u \|_{\lambda,\sigma,s}^2 &+ \| u \|_{\lambda,\sigma+1,s}^2 \lesssim \left(\| \na_x u\|_{\infty}+\| f\|_{\sigma,M}+\lambda^2 \| f\|_{\lambda,\sigma,M,s}\right) \| u \|_{\lambda,\sigma,s}^2 + \| u \|_\sigma^2\| u \|_{\lambda,\sigma,s}  \nonumber \\
 &+\| f\|_{\lambda,\sigma,M,s} \| u \|_{\lambda,\sigma,s}+ \left( \dot{\lambda}+\lambda \| u\|_\sigma + \lambda^2 (\| u \|_\sigma+\|u\|_{\lambda,\sigma,s}) \right) \| u \|_{\lambda,\sigma+\frac{s}{2},s}^2.
\end{align*}
In particular, by choosing $\lambda$ such that
\begin{equation}\label{condition lambda}
\dot{\lambda}+\lambda (1+\|f\|_{\sigma,M}+\| u\|_\sigma) +\lambda^2 (\|u\|_\sigma+\| f \|_{\lambda,\sigma,M,s}+\| u \|_{\lambda,\sigma,s}) \leqslant 0 ,
\end{equation}
we obtain:
\begin{equation}\label{estimation f sans lambda}
\frac{1}{2}\frac{\ud}{\ud t}  \| f \|_{\lambda,\sigma,M,s}^2 \lesssim \|f\|_{\lambda,\sigma,M,s} \bigg( \left(\|u\|_{W^{1,\infty}}+1\right) \|f\|_{\lambda,\sigma,M,s} + \|u\|_{\sigma} \|f\|_{\sigma,M} \bigg)
\end{equation}
and 
\begin{equation}\label{estimation u sans lambda}
\frac{1}{2}\frac{\ud}{\ud t} \| u \|_{\lambda,\sigma,s}^2 \lesssim  \| u \|_{\lambda,\sigma,s}\bigg( \left(\|\na_x u\|_{\infty}+\| f\|_{\sigma,M}\right) \| u \|_{\lambda,\sigma,s} + \| u \|_\sigma^2 + \| f\|_{\lambda,\sigma,M,s} \bigg)  . 
\end{equation}
Note that we will look for $\lambda$ such that 
\begin{equation}\label{condition lambda = 0}
\dot{\lambda}+\lambda (1+\|f\|_{\sigma,M}+\| u\|_\sigma) +\lambda^2 (\|u\|_\sigma+\| f \|_{\lambda,\sigma,M,s}+\| u \|_{\lambda,\sigma,s}) = 0.
\end{equation}
The last two inequalities \eqref{estimation f sans lambda} and \eqref{estimation u sans lambda} lead to the following two
\begin{equation}\label{G-V1}
\frac{\ud}{\ud t}  \| f \|_{\lambda,\sigma,M,s} \leqslant  C_1'  \left(\|u\|_{W^{1,\infty}}+1\right) \|f\|_{\lambda,\sigma,M,s} + C_1' \|u\|_{\sigma} \|f\|_{\sigma,M} 
\end{equation}
and 
\begin{equation}\label{G-NS1}
\frac{\ud}{\ud t} \| u \|_{\lambda,\sigma,s} \leqslant C_2'  \left(\|\na_x u\|_{\infty}+\| f\|_{\sigma,M}\right) \| u \|_{\lambda,\sigma,s} + C_2'\left(\| u \|_\sigma^2 + \| f\|_{\lambda,\sigma,M,s}\right) , 
\end{equation}
We will first deal with inequality \eqref{G-V1}. Once we have a bound for $\| f\|_{\lambda,\sigma,M,s}$, we integrate inequality \eqref{G-NS1}.  By Proposition \ref{Sobolev VNS}, we have: 
$\ \ds \|f\|_{\sigma,M}^2+\|u\|_\sigma^2 \leqslant \left(\|f_0\|_{\sigma,M}^2+\|u_0\|_\sigma^2\right) g(t),$
with
$$ g(t):=  \exp\big[C_0 \int_0^t \left(\|u(\tau)\|_{W^{1,\infty}}+\| \rho(\tau)\|_\infty+\| f(\tau) \|_{\infty,M}^2+1 \right)\ud \tau\big].$$
So if for $T>0$ the quantity $g(T)$ is finite then,  for all $t \in [0,T]$
$$ \int_0^t \|u(\tau)\|_{\sigma} \|f(\tau)\|_{\sigma,M} \ud \tau < \infty \quad \mbox{ and } \quad\int_0^t \left(\|u(\tau)\|_{W^{1,\infty}}+1\right) \ud \tau < \infty . $$ 
Therefore, \eqref{G-V1} with Gronwall's inequality lead us to
\begin{equation}\label{G for f}
 \| f(t) \|_{\lambda,\sigma,M,s} \leqslant \left( \|f_0\|_{\lambda_0,\sigma,M,s} + C_1' \int_0^t \|u(\tau)\|_{\sigma} \|f(\tau)\|_{\sigma,M} \ud \tau \right) \ue^{ C_1' \int_0^t \left(\|u(\tau)\|_{W^{1,\infty}}+1\right) \ud \tau } .  
\end{equation}
Note that the previous inequality together with $\|u\|_{W^{1,\infty}} \lesssim \|u\|_{\sigma}$ and Proposition \ref{Sobolev VNS} imply that  
\begin{equation}\label{G-V}
\| f \|_{\lambda,\sigma,M,s} \leqslant   C_1 \left(t+1\right) g(t)  .
\end{equation}
Returning to inequality \eqref{G-NS1},  since we already have $\sigma > \frac{d}{2}+1$ and $\| u \|_{\sigma} \leqslant \| u \|_{\lambda,\sigma,s}$, we write:
\begin{align*}
\frac{\ud}{\ud t} \| u \|_{\lambda,\sigma,s} &\lesssim \left(\|\na_x u\|_{\infty}+\| f\|_{\sigma,M}\right) \| u \|_{\lambda,\sigma,s} + \left(\| u \|_\sigma^2 + \| f\|_{\lambda,\sigma,M,s}\right) \\
&\lesssim \left(\|u\|_{\sigma}+\| f\|_{\sigma,M}\right) \| u \|_{\lambda,\sigma,s}  + \| f\|_{\lambda,\sigma,M,s} .
\end{align*}
 Hence,  by Gronwall's inequality we obtain:
\begin{equation}\label{G for u}
 \|u(t)\|_{\lambda,\sigma,s} \leqslant \left(\|u_0\|_{\lambda_0,\sigma,s} + C \int_0^t \|f(\tau)\|_{\lambda,\sigma,M,s}  \ud \tau \right)\ue^{C\int_0^t  \left(\|u(\tau)\|_{\sigma}+\| f(\tau)\|_{\sigma,M}\right)  \ud \tau}.
\end{equation}
Thanks to inequality \eqref{G-V} and Proposition \ref{Sobolev VNS}, we can write:
\begin{equation}\label{G-NS}
\| u \|_{\lambda,\sigma,s} \leqslant  \left(\| u_0 \|_{\lambda_0,\sigma,s}+C_2\int_0^t (1+\tau) g(\tau)\ \ud \tau \right)\ue^{C_2 \int_0^t g(\tau) \ud \tau } ,
\end{equation}
where $C_1$ and $C_2$ in \eqref{G-V} and \eqref{G-NS} respectively,  are two positive constants which depend on the initial data $(f_0,u_0)$, the radius of regularity $\lambda_0$, the Sobolev correction $\sigma$, the weight $M$ and the dimension $d$.  

\begin{remark}
\item \begin{enumerate}
\item The estimates \eqref{G for f} and \eqref{G for u} are valid for $s=1$.  
\item Note that we can estimate $ \| f \|_{\lambda,\sigma,M,s}^2 +\| u \|_{\lambda,\sigma,s}^2$ directly, by summing the two inequalities \eqref{estimation f sans lambda} and \eqref{estimation u sans lambda}, using $\|u\|_ {W^{1,\infty}} \lesssim \|u\|_{\sigma}$ for $\sigma > \frac{d}{2}+1$, $\|u\|_{\sigma } \leqslant \| u\|_{\lambda,\sigma,s}$ and $\| f\|_{\sigma,M} \leqslant \| f\|_{\lambda,\sigma,M,s}$, in order to get:
\begin{equation}\label{estimation Gevrey de f + u}
\| f \|_{\lambda,\sigma,M,s}^2 +\| u \|_{\lambda,\sigma,s}^2 \leqslant \big(\| f_0 \|_{\lambda_0,\sigma,M,s}^2 +\| u_0 \|_{\lambda_0,\sigma,s}^2\big) \ \exp\big[ \tilde C_1 \int_0^t \big(1+\| f\|_{\sigma,M}+\|u\|_\sigma\big) \ud \tau\big] .
\end{equation}
\end{enumerate}
\end{remark}
Finally, we need to find a positive continuous function $\lambda$ such that inequality \eqref{condition lambda} holds. 
\begin{lemma}\label{lemme lambda}
For Gevrey initial data $(f_0,u_0)$, with $\lambda_0 > 0$, $\sigma > \frac{d}{2}+\frac{s}{2}+2$ and $M > \frac{d}{2}+1$,  the function $\lambda$ defined on $[0,T_{max})$ by
\begin{equation}\label{lambda}
\lambda(t) := \frac{1}{G(t)}\left( \lambda_0^{-1} + \int_0^t \frac{\|u\|_\sigma+\| f \|_{\lambda,\sigma,M,s}+\| u \|_{\lambda,\sigma,s}}{G(\tau)} \ud \tau \right)^{-1} 
\end{equation}
satisfies the condition \eqref{condition lambda}, and where the function $G$ in the previous formula is given by
\begin{equation}\label{G(t)}
G(t) := \exp \left[ \int_0^t \left( 1+ \|u\|_\sigma+\| f \|_{\sigma,M} \right) \ud \tau \right] .
\end{equation}
Moreover, $\lambda$ is positive and non-increasing and satisfying
\begin{equation}\label{lambda(t)>...}
\lambda(t) \geqslant  \left(2C_3t+\lambda_0^{-1} \right)^{-1} \exp \bigg[-C_3 \int_0^t \left( 1+ \|u\|_\sigma+\| f \|_{\sigma,M} \right) \ud \tau \bigg] > 0 ,
\end{equation}
where $C_3$ is a constant depending on the initial data $(f_0,u_0)$ , the radius of regularity $\lambda_0$, the Sobolev correction $\sigma$, the weight $M$ and the dimension $d$. 
\end{lemma}

\ni\bpl \ref{lemme lambda}. The function $\lambda$ given by the formula \eqref{lambda} is solution to the differential equation \eqref{condition lambda = 0}. Also, thanks to inequality \eqref{estimation Gevrey de f + u}
$$ \| f \|_{\lambda,\sigma,M,s} +\| u \|_{\lambda,\sigma,s} \leqslant \tilde C_2 \ \exp\big[ \tilde C_2 \int_0^t \big(1+\| f\|_{\sigma,M}+\|u\|_\sigma\big) \ud \tau\big] ,  $$ 
and since the function $G$, given by \eqref{G(t)} is nondecreasing then, from \eqref{lambda} we get:
\begin{align*}
\lambda(t) &=  \left( \lambda_0^{-1}G(t) + \int_0^t \big(\|u\|_\sigma+\| f \|_{\lambda,\sigma,M,s}+\| u \|_{\lambda,\sigma,s}\big)\frac{G(t)}{G(\tau)} \ud \tau \right)^{-1} \\
&\geqslant \left( \lambda_0^{-1}G(t) + 2 \tilde C_2 \int_0^t G^{\tilde C_2}(\tau) \frac{G(t)}{G(\tau)} \ud \tau \right)^{-1}  \\
&\geqslant \left(  \lambda_0^{-1} + 2Ct \right)^{-1} G^{-C}(t) \\
&= \left( 2Ct+ \lambda_0^{-1}\right)^{-1} \ue^{-C \int_0^t (1+\| f\|_{\sigma,M}+\|u\|_\sigma ) \ud \tau}   ,
\end{align*} 
where $C := \max\{1;\tilde C_2\}$.
\epl
$\bullet$ For a Gevrey initial data $(f_0,u_0)$,  we have  in particular,  
$$ \|u_0\|_\sigma^2 + \| f_0 \|_{\sigma,M}^2 \leqslant \| f_0 \|_{\lambda_0,\sigma,M,s}^2+\| u_0 \|_{\lambda_0,\sigma,s}^2 < + \infty .$$
Then, by Proposition \ref{Sobolev VNS} and for $g(t) < +\infty$, 
$$ \| f(t) \|_{\sigma,M}^2 + \|u(t)\|_\sigma^2  \leqslant \big(\|u_0\|_\sigma^2 + \| f_0 \|_{\sigma,M}^2 \big)  g(t)  < + \infty .$$
Which implies that $\quad \ds \| f(t) \|_{\lambda,\sigma,M,s}^2 + \| u(t) \|_{\lambda,\sigma,s}^2 < + \infty .$
\ep
 
\ni \textbf{Proof of Theorem \ref{propagation Analyticite}.} Recall that the only difference between the cases $s<1$ and $s=1$ is in the Sobolev estimates, Lemma \ref{lem Sob V} requires the assumption $f$ to be compactly supported in velocity, unlike Corollary \ref{corollaire estimation Sobolev}. Thus, it remains to estimate the Sobolev norms without the compact support condition for $f$.
By using the Sobolev embeddings $\|u\|_{W^{1,\infty}} \lesssim \|u\|_{\sigma}$ and $\|\rho_f\|_\infty \lesssim \|\rho_f\|_\sigma \lesssim \|f\|_{\sigma,M}$ for $\sigma > \frac{d}{2}+1$ and $M>\frac{d}{2}$ then, the two inequalities \eqref{estimation Sobolev pour f} and \eqref{Sobolev-NS2} become
$$ \frac{1}{2}\frac{\ud}{\ud t}  \|f\|_{\sigma,M}^2  \lesssim \big(\|u\|_{\sigma}+1\big) \|f\|_{\sigma,M}^2 $$
and 
 $$
 \frac{1}{2}\frac{\ud}{\ud t}  \|u\|_{\sigma}^2 +\|u\|_{\sigma+1}^2 \lesssim \big(\|u\|_{\sigma}+\|f\|_{\sigma,M}+1\big) \big(\|u\|_{\sigma}^2 + \|f\|_{\sigma,M}^2 \big) 
 $$ 
respectively. Thus, by summing these last two inequalities and using Young's inequality, we get
$$  \frac{1}{2}\frac{\ud}{\ud t} \big( \|u\|_{\sigma}^2 + \|f\|_{\sigma,M}^2 \big) + \|u\|_{\sigma+1}^2 \leqslant C \big(1+\|u\|_{\sigma}^2 + \|f\|_{\sigma,M}^2\big)^{\frac{3}{2}} .  $$
Therefore, by integrating this last inequality with the notations $Z(t) := \|u(t)\|_{\sigma}^2 + \|f(t)\|_{\sigma,M}^2$ and $Z_0:=Z(0)$, we obtain 
\begin{equation}\label{estimation Sobolev s=1}
 Y(t) := \|u(t)\|_{\sigma} + \|f(t)\|_{\sigma,M} \leqslant \sqrt{2}\big(1 + Z(t) \big)^{\frac{1}{2}}  \leqslant \sqrt{2}\frac{(1+Z_0)^{\frac{1}{2}}}{1-Ct(1+Z_0)^{\frac{1}{2}}} 
 \end{equation}
 for $$ t < T_0 := T_0(u_0,f_0) := C^{-1}(1+Z_0)^{-\frac{1}{2}}  .  $$
 Now if we know that $Z(T_0) < \infty $ then we can repeat the argument above.   \\

\ni $\bullet$ For the propagation of analyticity, we proceed as in the Gevrey case. Indeed, by inequality \eqref{G-V1}  and since we already have $\sigma > \frac{d}{2}+1$, we write 
 $$ \frac{\ud}{\ud t}  \|f\|_{\lambda,\sigma,M,s}  \lesssim \big(\|u\|_{W^{1,\infty}}+1\big) \|f\|_{\lambda,\sigma,M,s} + \|u\|_\sigma \|f\|_{\sigma,M} \lesssim \big(\|u\|_\sigma +1 \big) \|f\|_{\lambda,\sigma,M,s} . $$
 Hence, 
$$
 \|f(t)\|_{\lambda,\sigma,M,s} \leqslant \|f_0\|_{\lambda_0,\sigma,M,s} \ \exp\big[ C\int_0^t (1+\|u(\tau)\|_\sigma) \ud \tau \big].
$$
For $u$,  inequality \eqref{G for u} is valide for $s=1$ and one has with the notation $Y(t) := \|u\|_\sigma + \|f\|_{\sigma,M}$
$$
 \|u(t)\|_{\lambda,\sigma,1} \leqslant \left(\|u_0\|_{\lambda_0,\sigma,1} + C \int_0^t \|f(\tau)\|_{\lambda,\sigma,M,s}  \ud \tau \right)\ue^{C\int_0^t Y(\tau) \ud \tau}. 
 $$
\ep

\ni \textbf{Proof of Corollary \ref{existence globale}.} By Theorem \ref{thm HMM}, for initial modulated energy $\mathcal{E}(0)$ small enough, in the sense of \eqref{petitesse de E}, we get
$$ \int_0^\infty \| \na_x u(\tau)\|_{L^\infty(\mathbb{T}^3)} \ud \tau + \| \rho_f \|_{L^\infty((0,\infty)\times\mathbb{T}^3)} < \infty .$$
Then, for initial data $(f_0,u_0) \in \mathcal{G}^{\lambda_0,\sigma,M,\frac{1}{s}}(\mathbb{T}^3\times\RR^3)\times\mathcal{G}^{\lambda_0,\sigma,\frac{1}{s}}(\mathbb{T}^3)$ such that $f_0$ has a compact support in velocity, the quantity
$$ \exp\left[C_1 \int_0^t \left(\|u(\tau)\|_{W^{1,\infty}}+\|\rho(\tau)\|_\infty+\| f(\tau) \|_{\infty,M}^2+1\right)\ud \tau\right]$$
is finite for every $t \geqslant 0$. This implies that the Sobolev norms, and consequently the Gevrey norms, are finite for all $t \geqslant 0$. Hence the end of the proof.
\ep

To complete this section, we give the following lemma on the growth of the norm $\|f(t)\|_{\infty,M}$. 
\begin{lemma}\label{lemma f_inf,M}
Let $t > 0$. If $\|f_0\|_{\infty,M}$ is finite and $u \in L^1(0,t;L^\infty(\mathbb{T}^d))$, then $f(t)$ satisfies 
\begin{equation}\label{estimation f_inf,M}
\big\|f(t)\big\|_{\infty,M} \lesssim \ds \ue^{d t} \big(1+\|u\|_{L^1(0,t;L^\infty(\mathbb{T}^d))}^M\big) \big\|f_0\big\|_{\infty,M} .
\end{equation}
Moreover, for $d=3$, one has under the assumptions of Theorem \ref{thm HMM} 
\begin{equation}
\big\|f(t)\big\|_{\infty,M} \lesssim \ds \ue^{3t} \big\|f_0\big\|_{\infty,M} .
\end{equation}
\end{lemma} 
\bp First of all, note that $$\|f(t)\|_{\infty,M} := \bigg(\underset{|\alpha| \leqslant M}{\sum}\|v^\alpha f\|_\infty^2\bigg)^{\frac{1}{2}} \lesssim \underset{x,v}{\sup}\big(1+|v|^M\big) f(t,x,v) \lesssim  \|f(t)\|_{\infty,M} . $$
The proof of this lemma follows from that of Lemma 4.6 in \cite{HMM} and from some estimates obtained in the same reference. Indeed, by integrating the differential equation \eqref{caracteristiques}$_2$ satisfied by $s \mapsto V(s,t,x,v)$, we get
$$ V(0,t,x,v) = v \ue^t - \int_0^t \ue^s u \big(s,X(s,t,x,v)\big) \ud s.  $$  
Therefore,
\begin{align}\label{v < V + int}
|v| &\leqslant |V(0,t,x,v)| \ue^{-t} + \int_0^t \ue^{s-t} \|u(s)\|_\infty \ud s \\
& \leqslant |V(0,t,x,v)| +\|u\|_{L^1(0,t;L^\infty(\mathbb{T}^d))} .  \nonumber
\end{align}
Hence inequality \eqref{estimation f_inf,M} holds thanks to the identity $ f(t,x,v) = \ue^{d t} f_0\big(X(0,t,x,v),V(0,t,x,v)\big)$:
$$\big(1+|v|^M\big) f(t,x,v) \lesssim  \underset{x,v}{\sup}\big(1+|v|^M\big) f_0(x,v) \lesssim  \|f_0\|_{\infty,M} .$$
For $d=3$,  from \eqref{v < V + int} and by Sobolev's embedding $H^2(\mathbb{T}^3)  \hookrightarrow  L^\infty(\mathbb{T}^3)$ and the Cauchy-Schwarz inequality, we write:
\begin{align*}
|v| &\lesssim |V(0,t,x,v)|  + \left(\int_0^t \ue^{s-t} \|u(s)\|_{H^2}^2 \ud s\right)^{\frac{1}{2}} \\
&\lesssim |V(0,t,x,v)|  + \underset{[0,t]}{\sup}\ \|u\|_{L^2} + \left(\int_0^t \ue^{s-t} \left(\big\|\na_x u(s)\big\|_{L^2}^2 + \big\|\Delta_x u(s)\big\|_{L^2}^2\right) \ud s\right)^{\frac{1}{2}} .
\end{align*}
Finally,   $ |v| \lesssim |V(0,t,x,v)| + 1$, thanks to the decay of the energy \eqref{decay E},   inequality (5.5) of \cite[Proposition 5.3]{HMM} namely:
$$ \big\|\na_x u(t)\big\|_{L^2}^2 + \int_{1/2}^t\big\|\Delta_x u(s)\big\|_{L^2}^2 \ud s \lesssim E(0) \big(1+\underset{[0,t]}{\sup} \ \|\rho_f\|_{\infty}\big) ,$$
and the smallness of $\mathcal{E}(0)$ which ensures that all the previous quantities are finite by Theorem \ref{thm HMM}. 
\ep

\ni \textbf{Comments.}
\begin{enumerate}
\item These last estimates on $\|f(t)\|_{\infty,M}$ imply that this norm grows at most as $\exp(t)$, which implies that the Sobolev norm grows at most as $\exp(\exp(t))$ and consequently, the Gevrey norm as a triple exponential. This is because of the term $\ue^{dt}$ which comes from $\na_v\cdot(vf)$. Without this term we would have found exactly the same estimates for Vlasov as in \cite{VR}.
\item We have managed to obtain estimates of the Sobolev norms over a time interval that does not depend on the data, thanks to the dissipative term in the NS equations and thanks to inequality 2 of Lemma \ref{commutateur}, which requires $f$ to be compactly supported in velocity, which is analogous to the result of \cite{VR}, where the force $F$ satisfies the Poisson equation. But the difficulty here is that the force depends on $u$, and $u$ itself still depends on $f$.  
\item Since we have products, $uf$ and $\rho u$, we are obliged to take the same radius of regularity/analyticity, which reduces the chances of taking advantage of the regularity of the NS solutions. 
\item The smallness condition given in \cite{HMM} is sufficient to control the quantities $\|\rho_f\|_{\infty}$ and $\int_0^\infty\|\na_x u(t)\|_\infty \ud t$. To the author's knowledge, it is not known whether this condition is optimal.  Similarly, the author does not know whether the exponential loss of the radius of regularity/analyticity is optimal, but what is certain is that $\lambda$ decays at least exponentially, because the Gevrey estimate that comes from the transport terms of the Vlasov equation is sharp. See estimates of $\hat E$ and $\hat F$ in Section \ref{sec Gevrey Vlasov}.
\end{enumerate}
 
 \appendix
 \section{Appendix}
 \subsection{Proof of Lemma \ref{inegalites triangulaires}}\label{appendice}
Some of inequalities of Lemma \ref{inegalites triangulaires} are found in \cite{VR} and some others in \cite{BMM}, and in the two cited references the proof is not given.   Let $k,l,\eta \in \RR^d$.\\
\textbf{1.} First, let show that $\ds \langle k+l\rangle \leqslant \langle k\rangle +\langle l\rangle$.  For this, we consider the increasing convex function $f:[0,+\infty[ \longrightarrow \RR$ defined by $$f(t):=(1+t^2)^{\frac{1}{2}} . $$
Observe that $\ds f(|k|) = \langle k\rangle$. We have
$$ \frac{1}{4} \langle k+l\rangle = \frac{1}{4} f(|k+l|) \leqslant \frac{1}{4} f(|k|+|l|) \leqslant \frac{1}{2} f\left(\frac{|k|+|l|}{2}\right) \leqslant \frac{f(|k|)+f(|l|)}{4} = \frac{1}{4} (\langle k\rangle + \langle l\rangle ) .$$ 
Now, let $s\in[0,1]$ and let consider the increasing function $g:[1,+\infty[ \longrightarrow \RR$ defined by $$ g(t):= 1+t^s-(1+t)^{s} . $$
We have for all $t \geqslant 1$, $g(t) \geqslant g(1)=2-2^s \geqslant 0$. Thus, $1+t^s\geqslant (1+t)^{s}$ for all $t \geqslant 1$. Let $k,l \in \RR^d$. Without loss of generality, we can assume that $|k| \geqslant |l|$. Therefore, for $t=\frac{\langle k\rangle}{\langle l\rangle} \geqslant 1$, we obtain
$$ \left( 1+ \frac{\langle k\rangle}{\langle l\rangle} \right)^s \leqslant 1 + \frac{\langle k\rangle^s}{\langle l\rangle^s}.$$
Which implies that
$$ \langle k+l \rangle^s \leqslant \left( \langle k\rangle + \langle l\rangle \right)^s \leqslant \langle k\rangle^s + \langle l\rangle^s, $$
since the function $t \mapsto t^s$ is increasing for $t \geqslant 1$ and $s\geqslant0$, and we have $\ds \langle k+l\rangle \leqslant \langle k\rangle +\langle l\rangle$.\\

\ni \textbf{2.} The function $t \mapsto t^s$ with $t \geqslant 1 $, is convex for all $s\geqslant 1$. Therefore, for $a\geqslant 1$ and $b \geqslant 1$
$$ \left( \frac{a+b}{2} \right)^s \leqslant \frac{a^s+b^s}{2}.$$
Thus, for $a=\langle k \rangle$ and $b=\langle l \rangle$, we get
$$ \langle k+l \rangle^s \leqslant \left( \langle k\rangle + \langle l\rangle \right)^s \leqslant 2^{s-1} \big( \langle k\rangle^s + \langle l\rangle^s \big) . $$
\textbf{3.} Let $s \in [0,1]$. We have
$$
\left\{\begin{array}{lcr}
\langle k \rangle^s = \langle k-l+l \rangle^s \leqslant \langle k-l \rangle^s + \langle l \rangle^s ,\\
\langle l \rangle^s = \langle k-l+l \rangle^s \leqslant \langle k-l \rangle^s + \langle k \rangle^s.
\end{array}\right.
$$
Hence, $$\big| \langle k \rangle^s - \langle l \rangle^s \big|\leqslant \langle k-l \rangle^s  .$$
\textbf{4.} Let $s \in [0,1]$. By applying the mean value theorem to the function $t \mapsto t^s$ between $X=\langle k \rangle$ and $Y=\langle l \rangle$, we obtain
$$\langle k \rangle^s - \langle l \rangle^s = s \big(\theta\langle k \rangle +(1-\theta) \langle l \rangle\big)^{s-1} \big(\langle k \rangle - \langle l \rangle\big) ,$$
where $\theta := \theta_{k,l} \in ]0,1[$. Thus, by concavity of the function $t \mapsto t^{s-1}$ for $s \in [0,1]$ and inequality of the previous point for $s=1$, we write
$$\big| \langle k \rangle^s - \langle l \rangle^s \big| \leqslant \frac{s\langle k-l \rangle}{\big(\theta\langle k \rangle +(1-\theta) \langle l \rangle\big)^{s-1}} \leqslant \frac{s\langle k-l \rangle}{\theta\langle k \rangle^{s-1} +(1-\theta) \langle l \rangle^{s-1}} \leqslant \frac{s}{\min(\theta,1-\theta)} \frac{\langle k-l \rangle}{\langle k \rangle^{s-1} + \langle l \rangle^{s-1}}.$$
\textbf{5.} Let $s \in [0,1]$. By applying the mean value theorem to the function $t \mapsto (t^{2/s}+|\eta|^2)^\frac{s}{2}$ between $X=\langle k \rangle^s$ and $Y=\langle k-l \rangle^s$, there exists $\theta := \theta_{k,l,\eta} \in ]0,1[$ such that
$$\big|\langle k,\eta \rangle^s - \langle k-l,\eta \rangle^s\big| = \bigg| \frac{\big(\theta\langle k \rangle^s +(1-\theta) \langle k-l \rangle^s\big)^{\frac{s}{2}-1}}{\big[|\eta|^2+(\theta\langle k \rangle^s +(1-\theta) \langle k-l \rangle^s)^{\frac{2}{s}}\big]^{1-\frac{s}{2}}} \big(\langle k \rangle^s - \langle k-l \rangle^s\big)\bigg|\leqslant \big| \langle k \rangle^s - \langle k-l \rangle^s \big| ,$$
since for $s\in[0,1]$,
$$ \bigg| \frac{\big(\theta\langle k \rangle^s +(1-\theta) \langle k-l \rangle^s\big)^{\frac{s}{2}-1}}{\big[|\eta|^2+(\theta\langle k \rangle^s +(1-\theta) \langle k-l \rangle^s)^{\frac{2}{s}}\big]^{1-\frac{s}{2}}}\bigg|\leqslant 1 .$$
Thus, inequality 5 follows from \textbf{3}.\\

\ni \textbf{6.} Let $s \geqslant 1$. We have by \textbf{2},
$$ \langle k,\eta \rangle^s = \big(1+|\eta|^2+|k|^2\big)^\frac{s}{2} \leqslant \big(1+|\eta|^2+2|k-l|^2+2|l|^2\big)^\frac{s}{2} \leqslant 2^{\frac{3s}{2}-1}\big(\langle k-l,\eta \rangle^s+\langle l \rangle^s\big) .$$
\textbf{7.} The proof of this point is identical to that of \textbf{4}, with $X=\langle k,\eta \rangle$ and $Y=\langle k-l,\eta \rangle$.  
\ep

\subsection{Energy estimates and sufficient condition for global existence}
Here we recall the definitions of energy and modulated energy used in \cite{choi2015global} and \cite{HMM}.
\begin{definition}\label{energie}
The kinetic energy of the VNS system \eqref{VNS} is given for $t \geqslant 0$ by
\begin{equation}
E(t) := \frac{1}{2} \int_{\mathbb{T}^d} |u(t,x)|^2 \ud x + \frac{1}{2} \int_{\mathbb{T}^d\times\RR^d} |v|^2 f(t,x,v) \ \ud v\ud x ,
\end{equation}
and the dissipation is defined as
\begin{equation}
D(t) := \int_{\mathbb{T}^d\times\RR^d} |u(t,x)-v|^2 f(t,x,v)\ \ud v\ud x + \int_{\mathbb{T}^d} |\na_x u (t,x)|^2 \ud x .
\end{equation}
Moreover, the two preceding quantities satisfy the identity 
\begin{equation}
\frac{\ud}{\ud t} E(t) + D(t) = 0 .
\end{equation}
Therefore, for an admissible initial data $(u_0,f_0)$ in the sense of \cite[Definition 1.2]{HMM},  the following energy estimate holds for almost all $t \geqslant s \geqslant 0$,
\begin{equation}\label{decay E}
E(t) + \int_s^t D(s) \ud s \leqslant E(s) .
\end{equation}
\end{definition}
\begin{definition}\label{energie modulee}
We define the modulated energy as
\begin{align}
\mathcal{E}(t) :=  \frac{1}{2} \int_{\mathbb{T}^d\times\RR^d} |v&-\langle j_f(t) \rangle|^2 f(t,x,v) \ \ud v\ud x   + \frac{1}{2} \int_{\mathbb{T}^d} |u(t,x)-\langle u(t) \rangle|^2 \ud x \\
&+ \frac{1}{4} |\langle u(t) \rangle - \langle j_f(t) \rangle|^2 ,
\end{align}
where $\langle u(t) \rangle := \int_{\mathbb{T}^d} u(t,x) \ud x$ and $\langle j_f(t) \rangle := \int_{\mathbb{T}^d} j_f(t,x) \ud x$.
\end{definition}
\ni Now let $N_q(f)$ and $M_\alpha f$ denote the following two quantities:
$$ N_q\big(f(t)\big) := \underset{x \in \mathbb{T}^3,  v \in \RR^3}{\sup} \left(1+|v|^q\right) f(t,x,v) \quad \mbox{ and } \quad M_\alpha f (t) := \int_{\mathbb{T}^d\times\RR^d} |v|^\alpha f(t,x,v) \ \ud v \ud x .$$
With the previous definitions and notations, we have the following result which guarantees the control of $\|\rho_f\|_\infty$ and $\int_0^\infty \|\na_x u\|_\infty \ud \tau$ under the smallness condition on $\mathcal{E}(0)$.
\begin{theorem} [\textbf{\cite[Theorem 2.2]{HMM}}] \label{thm HMM}
There exists $C_\star > 0$ and a nondecreasing onto function $\varphi : \RR_+ \to \RR_+$ such that the following holds. Let $(u_0,f_0)$ be an admissible initial condition such that $N_q(f_0)<+\infty$ for some $q>4$, $M_\alpha f_0 <+\infty$ for some $\alpha>3$ and $u_0 \in  H^{1/2}(\mathbb{T}^3)$. Then,  if
$$ \| u_0 \|_{\dot{H}^{1/2}(\mathbb{T}^3)} < \frac{1}{C_\star} , $$
and if the initial modulated energy $\mathcal{E}(0)$ is small enough, in the sense that
\begin{equation}\label{petitesse de E}
\varphi \left(N_q(f_0) + M_\alpha f_0 + E(0) +  \| u_0 \|_{H^{1/2}(\mathbb{T}^3)} + 1  \right) \mathcal{E}(0) < \min\left(1, \frac{1}{C_\star^2 -  \| u_0 \|_{\dot{H}^{1/2}(\mathbb{T}^3)}^2}\right) ,
\end{equation}
then,  there exists $\lambda, C'_\lambda >0$ such that for all $t \geqslant 0$,
\begin{equation}
\mathcal{E}(t) \leqslant \mathcal{E}(0) C_\lambda' \ue^{-\lambda t}.
\end{equation}
Furthermore,  we have the global bounds
\begin{equation}
\underset{t \geqslant 0}{\sup}\ \| \rho_f(t)\|_{L^\infty(\mathbb{T}^3)} < + \infty\quad \mbox{ and } \quad \int_1^{+\infty} \| \na_x u(\tau)\|_{L^\infty(\mathbb{T}^3)} \ud \tau < \eta\big(\mathcal{E}(0)\big) ,
\end{equation}
where $\eta$  is a continuous nonnegative function such that $\eta(0) = 0$.
\end{theorem}

{\footnotesize
\newcommand{\etalchar}[1]{$^{#1}$}

}
\bibliographystyle{alpha}

\end{document}